# CUQIpy: I. Computational uncertainty quantification for inverse problems in Python

**Nicolai A B Riis**[1,4] ⓘ, **Amal M A Alghamdi**[1] ⓘ,
**Felipe Uribe**[2,5] ⓘ, **Silja L Christensen**[1] ⓘ,
**Babak M Afkham**[1] ⓘ, **Per Christian Hansen**[1] ⓘ
and **Jakob S Jørgensen**[1,3,*] ⓘ

[1] Department of Applied Mathematics and Computer Science, Technical University of Denmark, Richard Petersens Plads, Building 324, 2800 Kongens Lyngby, Denmark
[2] School of Engineering Sciences, Lappeenranta-Lahti University of Technology (LUT), Yliopistonkatu 34, Lappeenranta 53850, Finland
[3] Department of Mathematics, The University of Manchester, Oxford Road, Alan Turing Building, Manchester M13 9PL, United Kingdom
[4] Copenhagen Imaging ApS, Herlev, Denmark

E-mail: jakj@dtu.dk

## Abstract

This paper introduces CUQIpy, a versatile open-source Python package for computational uncertainty quantification (UQ) in inverse problems, presented as Part I of a two-part series. CUQIpy employs a Bayesian framework, integrating prior knowledge with observed data to produce posterior probability distributions that characterize the uncertainty in computed solutions to inverse problems. The package offers a high-level modeling framework with concise syntax, allowing users to easily specify their inverse problems, prior information, and statistical assumptions. CUQIpy supports a range of efficient sampling strategies and is designed to handle large-scale problems. Notably, the automatic sampler selection feature analyzes the problem structure and chooses a suitable sampler without user intervention, streamlining the process.

---

[5] Part of the work by F U was done while employed at the Technical University of Denmark.
[*] Author to whom any correspondence should be addressed.

With a selection of probability distributions, test problems, computational methods, and visualization tools, CUQIpy serves as a powerful, flexible, and adaptable tool for UQ in a wide selection of inverse problems. Part II of the series focuses on the use of CUQIpy for UQ in inverse problems with partial differential equations.

## 1. Introduction

Inverse problems arise in various scientific and engineering fields such as materials science, medical imaging, and geoscience, [8, 11, 20, 60]. In these problems we compute hidden or un-observable features from indirect measurements, and the result is often uncertain due to noise in the data and inaccuracies in the models [22, 29]. We need to characterize and evaluate this uncertainty, such that we can make informed and safe decisions based on the computed results.

The field of *uncertainty quantification* (UQ) for inverse problems takes its basis in a Bayesian formulation. This approach builds upon theory and methods from Bayesian inference that provide a powerful and flexible framework for UQ. For inverse problems, it integrates prior knowledge with observed data and the computational model to yield posterior probability distributions that characterize the uncertainty in the computed solution [2, 9, 29, 52, 55].

### 1.1. Computational UQ and CUQIpy

The field of UQ for inverse problems is in a phase of rapid growth due to the development of new theory and methods; see, e.g. [4, 56]. The same is true for *computational UQ for inverse problems*, that focuses on the development of efficient computational methods for performing UQ of large-scale inverse problems. Several software packages have already been developed for forward and inverse UQ, including UQLab [35], SIPPI [24] and MUQ [40]. The latter takes a general approach, but many of these packages often target specific applications and have limited generality. In particular few if any UQ software packages handling large-scale imaging problems such as x-ray computed tomography (CT) and image deblurring appear to be available.

For this reason we have developed CUQIpy (pronounced "cookie pie"), an open-source Python package, currently in version 1.0.0, for computational UQ that targets a range of imaging-type inverse problems and includes a number of efficient computational methods. With CUQIpy, the user can specify a Bayesian inverse problem (or use one of the built-in ones) and then perform UQ computations using a number of methods that we provide. The name of the software is derived from the research project CUQI, computational uncertainty quantification for inverse problems, that funds the software development.

CUQIpy builds on an abstraction layer aimed at helping non-experts in Bayesian inference—and at the same time we give expert users flexibility and full control of the computational methods. A key ingredient of CUQIpy is a high-level modelling framework for working with inverse problems in the Bayesian setting. This framework includes a syntax that closely matches the underlying mathematics and statistics, thus enabling users to specify their inverse problem, *a priori* information, and other statistical information in a concise and intuitive way. CUQIpy's framework also offers experienced users access to the 'machine room' such that they can modify or define the underlying sampling strategies, analysis methods, etc.

Probabilistic programming refers to software that facilitates the specification and inference of probabilistic problems, typically using Bayesian inference. General-purpose probabilistic programming languages such as PyMC [51], STAN [10], Pyro [5], Turing.jl [15] have been very successful in enabling users to focus on their modelling while leaving the sampling details to the library. Enabled by automatic differentiation (AD) these languages utilize gradient information of the posterior distribution for efficient sampling via highly optimized versions of the so-called No-U-Turn sampler (NUTS) [27]. However, for large-scale inverse problems where computational efficiency is crucial, more specialized sampling strategies that take different types of problem structure into account are often required. CUQIpy addresses this issue by supporting different sampling strategies that exploit various types of problem structures such as linearity of operators or relations between distributions.

CUQIpy has been designed to support large-scale imaging-type inverse problems, such as CT. In this way, it allows users to incorporate and extend existing implementations of their inverse problem models and plotting tools. It also has an array-agnostic modeling framework that enables users to substitute classic array libraries, such as NumPy [25], with other libraries, such as PyTorch [41, 42], to benefit from AD or GPU acceleration.

CUQIpy is designed as a stand-alone Python library providing core functionality that allows the user to model and solve a variety of inverse problems. We include a selection of probability distributions, test problems, computational methods, and visualization tools. Additional functionality, e.g. via third-party libraries, is available through a *plugin interface*. Currently, three CUQIpy plugins have been released, proving tools for CT with Core Imaging Library CIL [28, 39], automated differentiation with PyTorch [42] and finite-element modelling through FEniCS [33], see section 3.1.

### 1.2. A motivating example

To demonstrate the basic usage of CUQIpy we consider a 2D deconvolution problem, which is a linear inverse problem that we write as $Ax = y$. Here, the matrix $A$ (which is assumed known) is the *forward model*, and in case of 2D deconvolution $A$ represents blurring by a point spread function. The vector $y$ represents a random variable for the blurred and noisy image, of which we are given a particular observed realization $y^{\text{obs}}$. From this forward model and data we are to infer the sharp image represented by the random variable $x$. We load a particular example forward model and observed data from the collection of test problems in CUQIpy:

```
A, y_obs, info = Deconvolution2D(dim=256,
                                 phantom="cookie").get_components()
```

The underlying sharp image and the observed blurred and noisy image are shown in figure 1. Other inverse problems can be represented in the same way and CUQIpy provides tools for users to specify their own forward model—see the example in section 4.

CUQIpy provides a modelling language to express what is known or assumed about the parameters in terms of probability distributions in a so-called *Bayesian Problem*; see section 2.1. Here we assume additive Gaussian white noise on the data $y$ with zero mean and an unknown precision $s$ (inverse variance) which is assumed to follow a Gamma distribution. For the image $x$ we choose an edge-preserving Laplace Markov random field (LMRF) prior, cf section 2.3, with an unknown scale parameter $d^{-1}$, where $d$ is also assumed to follow a Gamma distribution. The Bayesian Problem then takes the form

$$d \sim \text{Gamma}\left(1, 10^{-4}\right), \tag{1a}$$

$$s \sim \text{Gamma}\left(1, 10^{-4}\right), \tag{1b}$$

$$\boldsymbol{x} \sim \text{LMRF}\left(\boldsymbol{0}, d^{-1}\right), \tag{1c}$$

$$\boldsymbol{y} \sim \text{Gaussian}\left(\boldsymbol{A}\boldsymbol{x}, s^{-1}\boldsymbol{I}\right). \tag{1d}$$

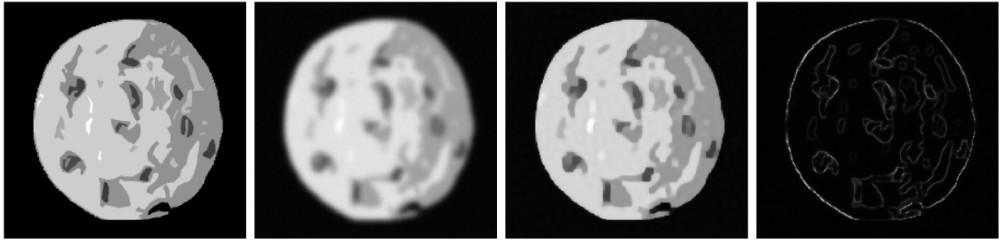

**Figure 1.** Figures generated by CUQIpy's automatic UQ analysis of the built-in 2D deconvolution test problem `Deconvolution2D` with a blurred cookie. Far left: true sharp image. Middle left: observed blurred and noisy image. Middle right: deblurred image (posterior mean). Far right: quantification of uncertainty (posterior width).

The parameters $d$ and $s$ are called hyperparameters (see section 2.1). In CUQIpy the syntax closely matches the mathematical specification of the Bayesian Problem:

```
d = Gamma(1, 1e-4)
s = Gamma(1, 1e-4)
x = LMRF(0, lambda d: 1/d, geometry=A.domain_geometry)
y = Gaussian(A @ x, lambda s: 1/s)
```

Here, lambda functions are used to define relations through algebraic operations on a parameter, e.g. $1/d$ and $1/s$; and the 'geometry' keyword to specify that the LMRF prior is defined on the domain of the forward operator A, which is represented by a so-called `Image2D` geometry—see section 3.

With the Bayesian Problem fully specified, we can immediately ask CUQIpy to perform a UQ analysis of the problem:

```
BP = BayesianProblem(d, s, x, y)   # Combine to Bayesian Problem
BP.set_data(y=y_obs)               # Specify observed data
BP.UQ()                            # Run UQ analysis
```

The `UQ()` method analyzes the problem, selects a suitable sampler, samples the posterior distribution, and returns a summary and selected visualizations of the results, as shown in figure 1. The deblurred image (posterior mean) demonstrates the impact of the LMRF prior in that it preserves the edge details in the image. The posterior width is also shown; providing an estimation of uncertainty in the deblurring, here highlighting that edges are subject to the most uncertainty.

In addition to estimating the image $\boldsymbol{x}$, also the hyperparameters $d$ and $s$ are estimated along with their uncertainty, which can be illustrated in various ways but is omitted here for brevity. The true value of $s$ is defined to be $7.716 \cdot 10^4$ in the `Deconvolution2D` test problem, and the posterior mean and 99% credibility interval are found to be $7.639 \cdot 10^4$ and $[7.535 \cdot 10^4, 7.744 \cdot 10^4]$. A true value for $d$ is not known but its posterior mean and credibility interval are found to be $3.871 \cdot 10^1$ and $[3.824 \cdot 10^1, 3.918 \cdot 10^1]$.

This illustrates the high-level usage of CUQIpy, where the user does not need to select a specific sampler or tune any parameters since CUQIpy does so automatically based on the problem structure. The following sections describe what happens 'under the hood' and how users can fully select and configure problem specification, samplers, visualization, etc.

### 1.3. Overview and notation

The present paper, which is the first of a two-part series, focuses on the core functionality of the CUQIpy package and its use for linear and nonlinear inverse problems. The companion paper [1] is dedicated to the modeling and UQ analysis of inverse problems based on partial differential equations (PDEs). We emphasize that both papers are written from a user's perspective, with a focus on descriptions of the software's functionality and examples of the use of CUQIpy.

Our paper is organized is follows. Section 2 briefly summarizes the theoretical background of Bayesian inference and computational UQ for inverse problems, and section 3 gives an overview of some important tools provided by CUQIpy demonstrated by a deconvolution problem. This is followed by three case studies that demonstrate the capabilities of CUQIpy: section 4 presents a case study of gravity anomaly detection with a nonlinear forward model illustrating how users can specify and solve their own inverse problem with CUQIpy. Section 5 describes a case study in x-ray CT using the CUQIpy-CIL plugin. Section 6 describes how the CUQIpy-PyTorch plugin expands CUQIpy with automatic differentiation from PyTorch and efficient sampling using Pyro and demonstrates this by specifying and solving the well-known benchmark problem 'Eight Schools' [16, 49]. Finally, section 7 discusses the current results and capabilities of CUQIpy, outlines future directions, and concludes the work.

We use the following notation: bold upper case such as $\boldsymbol{A}$ and $\boldsymbol{A}(\cdot)$ denote a linear and non-linear operator, respectively, with $\boldsymbol{I}_p$ denoting the $p \times p$ identity matrix; bold lower case such as $\boldsymbol{x}$ denotes a vector-valued random variable, and lower case such as $s$ and $f$ denotes a scalar random variable or a scalar function with $p$ denoting a probability density function; the superscript $k$ in $\boldsymbol{x}^{(k)}$ denotes the $k$th state in a Markov chain. A superscript text on a random variable such as $\boldsymbol{y}^{\text{obs}}$ denotes a realization.

## 2. Theoretical background

This section sets the stage for describing the software by summarizing important definitions and concepts in UQ; see [2, 13, 29] for more details and background.

### 2.1. Bayesian inverse problems

Many discretized inverse problems take their basis in the generic formulation

$$\boldsymbol{A}(\boldsymbol{x}) = \boldsymbol{y}, \tag{2}$$

where the vector $\boldsymbol{y}$ denotes the noisy data, the vector $\boldsymbol{x}$ represents the solution parameters to be found, and the operator $\boldsymbol{A}$ (which, in the linear case, is a matrix) represents the forward model. In the Bayesian setting, the solution parameters $\boldsymbol{x}$ and the data $\boldsymbol{y}$ are random variables. A statistical model for (2) is then characterized by the *joint probability distribution* $p(\boldsymbol{x},\boldsymbol{y})$ of the parameters and data.

A useful representation of the joint distribution is to list what is known or assumed about the parameters and data, as well as their relation through what we call a *Bayesian Problem*. In statistics this is sometimes referred to as a *generative model* [16, 17]. For the joint distribution $p(\boldsymbol{x},\boldsymbol{y})$ associated with (2), we define the generic Bayesian Problem in terms of two distributions:

$$\boldsymbol{x} \sim p(\boldsymbol{x}), \tag{3a}$$

$$\boldsymbol{y} \sim p(\boldsymbol{y}|\boldsymbol{x}). \tag{3b}$$

We use the terms *prior* and *data distribution* for the distributions associated with the solution parameters and the data, respectively, in (3). For brevity in our notation, the statistical dependence $\boldsymbol{y}|\boldsymbol{x}$ is omitted from the left-hand-side in (3b).

In Bayesian inverse problems, the goal is to infer the solution parameters given a particular realization of the data. The *posterior distribution* $p(\boldsymbol{x}|\boldsymbol{y})$ characterizes the distribution of the solutions $\boldsymbol{x}$ to the inverse problem, given the data $\boldsymbol{y}$. Via *Bayes' theorem* for continuous probability densities, we can express the posterior as

$$p(\boldsymbol{x}|\boldsymbol{y}) = \frac{p(\boldsymbol{x},\boldsymbol{y})}{p(\boldsymbol{y})} = \frac{p(\boldsymbol{y}|\boldsymbol{x})\,p(\boldsymbol{x})}{p(\boldsymbol{y})}. \tag{4}$$

Given fixed observed data $\boldsymbol{y}^{\text{obs}}$, $p(\boldsymbol{y}|\boldsymbol{x})$ considered as a function of $\boldsymbol{x}$ is known as the *likelihood function* or just *likelihood*, denoted $L(\boldsymbol{x}|\boldsymbol{y}=\boldsymbol{y}^{\text{obs}})$. Furthermore, $p(\boldsymbol{y})$ is a normalization constant that is usually omitted and we write the posterior as proportional to the product of the likelihood and prior:

$$p(\boldsymbol{x}|\boldsymbol{y}=\boldsymbol{y}^{\text{obs}}) \propto L(\boldsymbol{x}|\boldsymbol{y}=\boldsymbol{y}^{\text{obs}})\,p(\boldsymbol{x}). \tag{5}$$

In CUQIpy the user is generally only expected to (i) define statistical assumptions about parameters and data via a Bayesian Problem and (ii) provide observed data. CUQIpy will then automatically formulate the posterior distribution by combining likelihood(s) and prior(s) through Bayes' theorem.

## 2.2. Hyperparameters

It is often the case that the distributions depend on one or more unknown parameters. In the Bayesian paradigm, we can assign probability densities to them and include them in the Bayesian Problem. For example assuming the distributions in (3) depend on *hyperparameters* $d$ and $s$ respectively, the generic Bayesian Problem for the joint probability distribution $p(\boldsymbol{x},\boldsymbol{y},s,d)$ would be

$$d \sim p(d), \tag{6a}$$

$$s \sim p(s), \tag{6b}$$

$$\boldsymbol{x} \sim p(\boldsymbol{x}|d), \tag{6c}$$

$$\boldsymbol{y} \sim p(\boldsymbol{y}|\boldsymbol{x},s). \tag{6d}$$

The posterior associated with this Bayesian Problem becomes

$$p(\boldsymbol{x},d,s|\boldsymbol{y}=\boldsymbol{y}^{\text{obs}}) \propto L(\boldsymbol{x},s|\boldsymbol{y}=\boldsymbol{y}^{\text{obs}})\,p(\boldsymbol{x}|d)\,p(d)\,p(s), \tag{7}$$

meaning that the task is now to infer about $x$ as well as $d$ and $s$, which is indeed possible in CUQIpy. This is usually referred to as *hierarchical modeling* and we refer to such additional model parameters as *hyperparameters*. See, e.g. [2, section 5.2] for details. Another example in CUQIpy is given in section 3.5.

### 2.3. Markov random field priors

In addition to well-known distributions such as the Gaussian, gamma and inverse-gamma, we also need a few distributions associated with Markov random fields. In CUQIpy, these are specifically used for the modeling of certain priors. For definitions and details about Markov random fields, see [32].

Let $x_i$ denote the elements of the $n$-dimensional random vector $x$ and consider the differences $\triangle x_i \equiv x_i - x_{i-1}$ for $i = 1,\ldots,n$. Then, following [2, chapter 4], we say that $x$ is a zero-mean Gaussian Markov random field (GMRF) if $\triangle x_i \sim \text{Gaussian}(0, d^{-1})$, where $d$ is a precision parameter. With zero boundary conditions (other conditions are also possible), this implies that

$$x \sim \text{Gaussian}\left(\mathbf{0}, \mathbf{Q}^{-1}\right) \quad \text{with} \quad \mathbf{Q} = d\mathbf{L}_n, \quad \mathbf{L}_n = \begin{pmatrix} 2 & -1 & & & \\ -1 & 2 & -1 & & \\ & \ddots & \ddots & \ddots & \\ & & & -1 & 2 \end{pmatrix}, \quad (8)$$

which we for brevity denote $x \sim \text{GMRF}(\mathbf{0}, d)$. When used as a prior, this GMRF induces regularity (smoothness) on the samples—similar to first-order Tikhonov regularization [22, section 8.1]. Higher-order differences can also be used in a GMRF (see, e.g. [50]).

In some applications we prefer a prior that allows less regularity—giving samples that may resemble total-variation regularized solutions [22, section 8.6]. This is possible if we replace the Gaussian distribution for $\triangle x_i$ with a Laplace or Cauchy distribution, i.e. we assume $\triangle x_i \sim \text{Laplace}(0, b)$ or $\triangle x_i \sim \text{Cauchy}(0, b)$. Again assuming zero boundary conditions, the corresponding densities for $x$ are, respectively,

$$p(x) \propto \frac{1}{(2b)^{n+1}} \exp\left(-\frac{\|\mathbf{D}_n x\|_1}{b}\right) \quad \text{with} \quad \mathbf{D}_n = \begin{pmatrix} 1 & & & & \\ -1 & 1 & & & \\ & \ddots & \ddots & & \\ & & & -1 & 1 \\ & & & & -1 \end{pmatrix} \quad (9)$$

and

$$p(x) \propto \prod_{i=1}^{n+1} \frac{b}{b^2 + (x_i - x_{i-1})^2} \quad \text{with} \quad x_0 = x_{n+1} = 0. \quad (10)$$

When used as a prior $p(x)$ they are called Laplace and Cauchy MRF priors, respectively, and we use the short-hand notations $x \sim \text{LMRF}(\mathbf{0}, b)$ and $x \sim \text{CMRF}(\mathbf{0}, b)$, for a scalar parameter $b$.

If $X$ represents a square $N \times N$ image with pixel values $X_{ij}$ for $i,j = 1,\ldots,N$, and the vector $x$ of length $N^2$ consists of the columns of $X$ stacked on top of each other, then we can define a GMRF in which the horizontal and vertical differences between pixel values follow Gaussian$(0, d^{-1})$. In this case $n = N^2$ and, with zero boundary conditions, the precision matrix $\mathbf{Q}$ in (8) takes the form $\mathbf{Q} = d(\mathbf{I}_N \otimes \mathbf{L}_N + \mathbf{L}_N \otimes \mathbf{I}_N)$ in which $\otimes$ denotes the Kronecker product. Similarly, for LMRF and CMRF with zero boundary conditions the densities take the form

$$p(\boldsymbol{x}) \propto \frac{1}{(2b)^{n+1}} \exp\left(-b^{-1}\left(\tfrac{1}{2}\|(\boldsymbol{I}_N \otimes \boldsymbol{D}_N)\boldsymbol{x}\|_1 + \tfrac{1}{2}\|(\boldsymbol{D}_N \otimes \boldsymbol{I}_N)\boldsymbol{x}\|_1\right)\right) \quad (11)$$

and

$$p(\boldsymbol{x}) \propto \prod_{i=1}^{N+1}\prod_{j=1}^{N+1} \frac{b}{b^2 + (X_{ij} - X_{i-1,j})^2} \frac{b}{b^2 + (X_{ij} - X_{i,j-1})^2}, \quad (12)$$

with zero values in all boundary pixels $X_{ij}$ just outside the $N \times N$ domain. For more details, see [2, section 4.3] and [54].

### 2.4. Sampling methods in CUQIpy

Selecting the most suited sampling method for exploring the posterior in a particular problem—as well as configuring and running it to obtain high-quality samples—is not always straightforward. Indeed, one of the aims of CUQIpy is to aid the user in this by providing automated sampler selection and configuration for problems of recognized structure.

In a few cases a closed-form expression exists for the posterior, which can be used for direct sampling. A classic example is that of a Gaussian prior and Gaussian data distribution with linear forward model which results in the posterior being Gaussian as well. For problems of modest size one can exploit this for efficient direct sampling based on the analytical expression. This is automatically detected and used by CUQIpy in `BayesianProblem`.

In some cases it is possible to obtain a tractable and computationally efficient approximation of the posterior distribution. One such example is the *Laplace approximation* [19, 53], with which one obtains a Gaussian distribution centered at the maximum-a-posteriori (MAP) estimate having covariance of the inverse Hessian of the negative logarithm of the targeted posterior distribution.

In most cases, however, computation of samples are based on Markov chain Monte Carlo (MCMC) algorithms [44]; these produce a sequence of samples $\boldsymbol{x}^{(1)}, \boldsymbol{x}^{(2)}, \ldots$ that represents a *Markov chain*, whose stationary distribution is the target posterior. The initial part of the chain, called the burn-in, is discarded because these samples may not be representative of the desired distribution. Moreover, we typically thin the sequence of samples to reduce the correlation.

We now summarize some of the MCMC samplers that are available in CUQIpy, highlighting their unique features and which problems they can be used for. For reviews of the methods we refer to [7, 38, 44]. CUQIpy also provides several common MCMC diagnostics to help monitor convergence and independence of samples; these are described at the end of this section.

#### 2.4.1. Metropolis–Hastings (MH) sampling [26].
This classical method uses a two-stage procedure with proposal step and acceptance/rejection steps. The first step computes a proposal $\boldsymbol{x}'$ from the density $q(\boldsymbol{x}|\boldsymbol{x}')$ which is the conditional probability of $\boldsymbol{x}$ given the proposed state $\boldsymbol{x}'$. The second step computes the acceptance ratio

$$\alpha(\boldsymbol{x}, \boldsymbol{x}') = \min\left(1, \frac{p(\boldsymbol{x}')\, q(\boldsymbol{x}'|\boldsymbol{x})}{p(\boldsymbol{x})\, q(\boldsymbol{x}|\boldsymbol{x}')}\right), \quad (13)$$

which expresses the probability of accepting $\boldsymbol{x}'$. At state $k$ of the MH algorithm, the next state $\boldsymbol{x}^{(k+1)}$ is chosen by first sampling a candidate point $\boldsymbol{x}'$ from the proposal density $q(\cdot|\boldsymbol{x}^{(k)})$. Then $\boldsymbol{x}'$ is accepted with probability (13) and $\boldsymbol{x}^{(k+1)} = \boldsymbol{x}'$, otherwise it is rejected and $\boldsymbol{x}^{(k+1)} = \boldsymbol{x}^{(i)}$.

*2.4.2. Preconditioned Crank–Nicolson (pCN) sampling [12].* This method assumes that the prior has a Gaussian distribution. It uses the Crank–Nicolson finite-difference scheme to solve an underlying stochastic differential equation (SDE) that is invariant with respect to the posterior [21, 47]. When we choose the prior covariance matrix $S$ as preconditioner, we obtain the following mechanism for producing the proposed state $x'$ in the MH method, see [12]:

$$x' = \sqrt{1-s^2}\,x + s\,\xi, \quad \text{where} \quad \xi \sim \text{Gaussian}(\mathbf{0}, S), \quad s \in (0,1]. \tag{14}$$

Given a current state $x^{(k)}$, it follows from (14) that the associated proposal distribution is Gaussian with mean vector $\sqrt{1-s^2}\,x^{(k)}$ and covariance matrix $s^2 S$, where $s$ is a tunable parameter controlling the acceptance rate. We remark that while pCN is designed for a Gaussian prior, it can be extended to more general priors by applying a transformation to a standard Gaussian distribution.

*2.4.3. Gradient-based sampling methods based on the Langevin SDE [37, 46].* As mentioned above, the pCN method is based on the numerical solution of a certain SDE. There are, in fact, many sampling methods that take such an approach; here we focus on methods based on the overdamped *Langevin SDE*

$$\mathrm{d}x(t) = -\nabla f(x(t))\,\mathrm{d}t + \sqrt{2}\,\mathrm{d}w(t), \tag{15}$$

where $x \in \mathbb{R}^n$ and $w(t)$ is standard Brownian motion in $\mathbb{R}^n$. Under certain conditions, as $t \to \infty$ the solution $x(t)$ follows a distribution with density $p(x) \propto \exp(-f(x))$.

The connection between (15) and posterior sampling methods arises when we choose

$$f(x) \propto -\log p(x|y) \tag{16}$$

and apply a numerical scheme to solve the SDE, thus producing a sequence of discretized solutions $x^{(i)}$. The standard choice is the forward Euler methods which leads to a scheme known as the *unadjusted Langevin algorithm (ULA)*:

$$x^{(k+1)} = x^{(k)} - h\nabla f\!\left(x^{(k)}\right) + \sqrt{2h}\,z^{(k)}, \qquad z^{(k)} \sim \text{Gaussian}(\mathbf{0}, I_n), \tag{17}$$

where a new $z^{(k)}$ is drawn in each step, and $h$ is a step size that controls the convergence and discretization error. We note that ULA can be interpreted as using a proposal distribution $\text{Gaussian}(x^{(k)} - h\nabla f(x^{(k)}), 2hI_n)$ but without an acceptance/rejection step, hence the name *unadjusted* and we stress that it is not an MH method.

A related method, known as the *Metropolis-adjusted Langevin algorithm (MALA)*, applies the MH acceptance/rejection step (13) to the above-mentioned ULA proposal which guarantees that samples follow the posterior distribution. More advanced gradient-based methods are the Hamiltonian Monte Carlo method and its adaptive version the NUTS; we refer to [7, 34].

*2.4.4. Randomize-then-optimize (RTO) sampling [2, 3].* This method, which was originally designed for nonlinear problems, uses an acceptance/rejection step (13) similar to the MH method, but the proposal strategy is different in that it involves solution of a (nonlinear) least-squares problem.

Our current implementation targets linear problems (Gaussian posteriors) and does not require an acceptance/rejection step; we refer to this version as *linear RTO*. We draw an independent sample of the proposal $x'$ by solving a linear least-squares problem:

$$x' = \arg\min_{x} \|Mx - b\|_2^2 \quad \text{with} \quad M = \begin{bmatrix} A \\ C \end{bmatrix}, \quad b = \begin{bmatrix} y \\ C\mu \end{bmatrix} + \widetilde{z}, \tag{18}$$

where $\widetilde{z} \sim \text{Gaussian}(\mathbf{0}, I_{m+n})$. Moreover, $\mu$ is the prior mean and $C$ is a matrix such that $C^\top C$ equals the precision matrix $Q$ for the prior; typically, $C$ is the Cholesky factor of $Q$. The independence of the computed posterior samples depends on the accuracy of the algorithm used to solve (18); in CUQIpy we use the well-known CGLS iterative algorithm (see, e.g. [22, section 6.3.2]) which is suited for large-scale problems. The linear RTO algorithm can be used to compute samples from the *unadjusted Laplace approximation* method proposed in [58], which is tailored to linear problems with LMRF priors.

*2.4.5. Gibbs sampling [18].* This method is useful when the posterior is expressed as a joint distribution via conditional densities, and in connection with hierarchical models such as (7). Specifically, Gibbs sampling is useful when computing samples from the joint distribution is impractical, but drawing samples from the conditional distributions of given parameter components is feasible.

To illustrate Gibbs sampling, consider a generic joint distribution $p(x,y,z)$; this could e.g. be the joint posterior in (7). The next state in the chain for $x, y, z$ is generated from the previous state as follows:

$$x^{(k+1)} \sim p\left(x|y^{(k)}, z^{(k)}\right), \tag{19a}$$

$$y^{(k+1)} \sim p\left(y|x^{(k+1)}, z^{(k)}\right), \tag{19b}$$

$$z^{(k+1)} \sim p\left(z|x^{(k+1)}, y^{(k+1)}\right). \tag{19c}$$

The special case of sampling a single random vector involves treatment of the elements component-by-component, and this scheme is known as the *component-wise Metropolis–Hastings (CWMH)* algorithm. There are different strategies to choose the scanning order of a Gibbs sampler; see [14, 38].

*2.4.6. Conjugacy-based samplers.* Application of the Gibbs sampler is useful when we can easily sample from the posterior distribution of each parameter conditioned on all the other ones. For example, in case of a conjugate prior (where the posterior is in the same probability distribution family as the prior), the conditionals are often available in closed form. In some cases pairs of distributions can also be approximated by a conjugate relation; this is the case, e.g. for the Gaussian-LMRF pair [58]. In other cases it is necessary to use MCMC methods to sample (some of) the conditionals; in this case the Gibbs sampler is referred to as *hybrid*. A classical example is the so-called *Metropolis-within-Gibbs* algorithm [45]. One of the strengths of CUQIpy is to automatically detect when conditionals are easy to sample and use this within Gibbs sampling.

### 2.5. Common MCMC diagnostics

Samples computed via MCMC methods are dependent. The effect of dependencies on the accuracy of MC estimates can be quantified in terms of the autocorrelation between the samples, $\{x^{(i)}\}_{i=1}^{N_{\text{samp}}}$, assuming they have already reached stationarity. For dependent samples, the variance of an MC estimate—compared to estimates based on independent samples—is modified by a factor $\tau$ called the *integrated autocorrelation time* (IACT) given by

$$\tau = 1 + 2\sum_{\ell=1}^{N_{\text{samp}}} \rho(\ell), \tag{20}$$

where $N_{\text{samp}}$ is the number of samples, $\rho(\ell)$ is the sample estimate of the normalized autocorrelation function for lag $\ell$. Related to this is the *effective sample size* (ESS) which expresses the amount by which autocorrelation within the chain influences uncertainty in the computed estimates:

$$\text{ESS} = \frac{N_{\text{samp}}}{\tau}. \tag{21}$$

ESS and the so-called *R-hat diagnostic* [59] are available in CUQIpy provided via ArviZ [30]. For more background on MCMC diagnostics, we refer the reader to [48].

## 3. Specifying and solving Bayesian inverse problems with CUQIpy

In the following sections we describe the most important software components of CUQIpy. We first provide an overview of the package and then take a step-by-step approach through some fundamental tools to illustrate a typical process of specifying and solving a Bayesian inverse problem.

### 3.1. Overview

CUQIpy is a Python package for Computational UQ for Inverse Problems developed at the Technical University of Denmark (DTU). It is available from:

https://cuqi-dtu.github.io/CUQIpy/

along with information on how to install and get started, full documentation, and numerous demos and tutorials. The package is released under a permissive Apache v2.0 license and is developed fully open source on GitHub.

CUQIpy is object-oriented and consists of a number of classes representing the building blocks needed to specify and solve Bayesian inverse problems. This includes distributions, forward models, samplers and so-called *geometry* for representing the context of a problem e.g. whether 1D or 2D (figure 2). In addition, a number of plug-ins are available as separate packages that expand the functionality of CUQIpy:

- CUQIpy-CIL: a plugin for the Python package CIL [28, 39] providing access to forward models for x-ray CT. As of now 2D parallel-beam and fan-beam models are exposed, while 3D parallel-beam and cone-beam can be exposed from CIL when needed.
- CUQIpy-FEniCS: a plugin providing access to the finite element modelling language FEniCS [33], which is used here for solution of PDE-based inverse problems. This is illustrated in the companion paper [1].
- CUQIpy-PyTorch: a plugin providing access to the AD framework of PyTorch [42] within CUQIpy. It allows gradient-based sampling methods without manually providing derivative information of distributions and forward models. This plugin is illustrated in section 6.

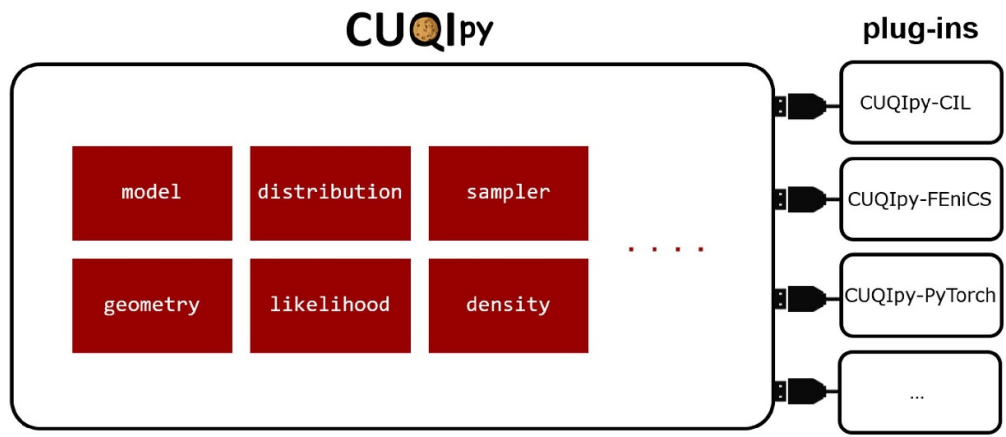

**Figure 2.** Diagram of CUQIpy with some of the main modules shown, as well as the plug-ins that each expand the core functionality.

**Table 1.** First block: test problems in CUQIpy not involving a PDE. Second block: test problems in CUQIpy involving a PDE and covered in companion article [1]. Third block: test problem in CUQIpy-CIL.

| Test problem name | Description |
| --- | --- |
| Deconvolution1D | 1D signal deblurring |
| Deconvolution2D | 2D image deblurring |
| Abel1D | Rotationally symmetric computed tomography |
| WangCubic | Problem with nonlinear two-parameter forward model |
| Heat1D | Discrete heat problem (time-dependent linear PDE) |
| Poisson1D | Discrete 1D Poisson problem (steady-state linear PDE) |
| ParallelBeam2D | 2D parallel-beam CT using CIL |

Once CUQIpy is installed, one may import the required components, for example—and for completeness—these are the classes needed for the examples in the present and introduction sections:

```
from cuqi.testproblem import Deconvolution1D, Deconvolution2D
from cuqi.distribution import Gaussian, Gamma, GMRF, LMRF, CMRF,
                              JointDistribution
from cuqi.problem import BayesianProblem
from cuqi.sampler import LinearRTO, UGLA, NUTS
import numpy as np
```

### 3.2. Test problems and forward models

CUQIpy provides a number of `testproblem` (cf table 1) which contain pre-made, configurable test problems. Those considered in this paper take the generic form $y = A(x)$, see (2), and they provide a forward model $A$, the underlying true signal $x^{\text{true}}$, as well as a clean data set $y^{\text{true}}$ and an observed noisy data set $y^{\text{obs}}$. Similar test problems, discussed in the companion paper [1], take their basis in a PDE formulation.

As running examples for this section we set up two `Deconvolution1D` linear test problems $y = Ax$ (recall also the `Deconvolution2D` test problem from section 1.2). First, we consider a smooth signal given by a sinc function which is subject to the default Gaussian blur, i.e. it is convolved with a Gaussian function:

```
A, y_obs, info = Deconvolution1D(phantom="sinc").get_components()
```

Second, we consider a piecewise constant signal (see documentation for further configuration options for signal type, blurring point spread function, noise level etc):

```
A, y_obs, info = Deconvolution1D(phantom="square",
                                 PSF_param=5).get_components()
```

Here we simply return the main components of the test problem, namely the forward model `A`, the observed data `y_obs` and additional information about the test problem in `info`. The `info` contains information about the test problem including the underlying true signal and exact data without noise that we extract here for convenience:

```
x_true = info.exactSolution
y_true = info.exactData
```

In CUQIpy, deterministic vectors such as $x^{\text{true}}$ or $y^{\text{true}}$ are represented by our `CUQIarray` data structure, which includes information about the type of signal in an attribute `geometry`. We can take a look at this for example for the clean signal

```
x_true.geometry
```

which for the 1D cases tells us that the signal has 128 elements and is of a 1D type:

```
Continuous1D(128,)
```

`CUQIarray` supports plotting, and the geometry also makes context-aware plotting simple by automatically making the correct type of plot associated with the class (see figures 1 and 3):

```
x_true.plot()
y_true.plot()
y_obs.plot()
(y_obs - y_true).plot()
```

In the last line we make use of the fact that `CUQIarray` is subclassed from conventional NumPy arrays meaning that all algebra is available; here we use that to determine and plot the noise that was added by the test problem to the clean blurred signal.

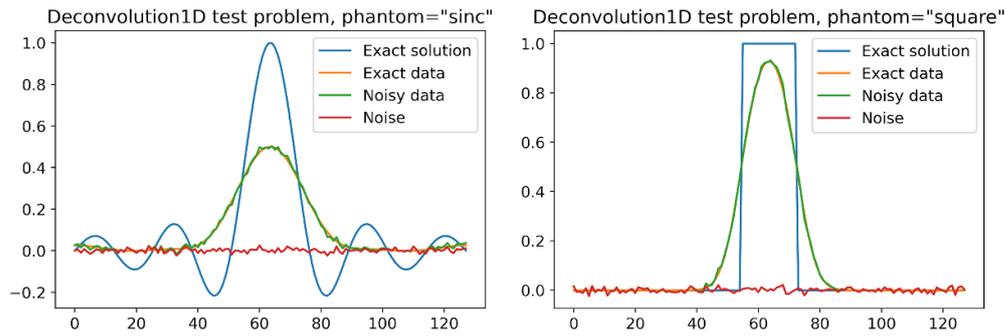

**Figure 3.** `Deconvolution1D` test problem with Gaussian blur. Left: a sinc function for the exact signal, the exact blurred and the blurred and noisy observed signal, and the additive noise. Right: Same for exact signal 'square'.

**Table 2.** Forward models in CUQIpy (top) and CUQIpy-CIL (bottom).

| Model name | Model description |
| --- | --- |
| `LinearModel` | Linear model $y = Ax$ |
| `Model` | General linear or non-linear model $y = A(x)$ |
| `PDEModel` | Model represented by a partial differential equation |
| `CILModel` | General model for CT problem with CIL |
| `ParallelBeam2DModel` | 2D parallel-beam scan model |
| `FanBeam2DModel` | 2D fan-beam scan model |
| `ShiftedFanBeam2DModel` | 2D fan-beam scan with shifted source/detector model |

The forward model `A` is of the type `LinearModel` which has methods `A.forward` and `A.adjoint` to compute forward and adjoint operations on a vector; as well as short hand NumPy-style matrix-vector multiplication syntax:

```
y_true = A @ x_true
y_adj  = A.T @ y_true
```

The last line explicitly creates and applies the adjoint operator, here the transposed matrix $A^\mathsf{T}$. CUQIpy handles non-linear forward models using the short hand syntax:

```
y_true = A(x_true)
```

Table 2 gives an overview of the types of models provided by CUQIpy as well as the CUQIpy-CIL plugin. Users can easily specify their own linear or nonlinear forward model representing the inverse problem of interest; an example is given in section 4.

Models are 'aware' of the spaces associated with the operator; in terms of geometries for the range and domain; these can be queried by

**Table 3.** Overview of common `geometry` classes in CUQIpy. A `geometry` represents the type of parameter, i.e. whether its elements are pixel values of an image or coefficients in an expansion etc and enables dedicated plotting according to signal type. The geometries are also responsible for transforming the underlying parameter vector to the so-called *function values* that are then passed to the forward model operator making the forward `Model` class agnostic to the array implementation.

| Geometry name  | Description                                         |
|----------------|-----------------------------------------------------|
| `Continuous1D` | Signal in one (space, time, etc) dimension          |
| `Continuous2D` | Signal in two (space or similar) dimensions         |
| `Image2D`      | Pixelated image signal                              |
| `Discrete`     | Collection of individual labelled scalar parameters |
| `MappedGeometry` | Maps another geometry by a callable function      |
| `KLExpansion`  | Karhunen–Loéve expansion of 1D signal               |
| `StepExpansion` | Piecewise constant/step expansion of 1D signal     |

```
A.domain_geometry
A.range_geometry
```

and since for the current deblurring examples the range and domain are identical both of these return the same geometry, in the 1D case `Continuous1D(128,)`. We can also query the dimensions and store these for later convenience, along with a zero vector:

```
m = A.range_dim
n = A.domain_dim
x_zero = np.zeros(n)
```

A variety of different geometries is provided (table 3) handling for example discrete labelled data as well as providing a way to parametrize solutions in terms of basis functions or a finite element mesh. Several examples of the use of geometries are given later (e.g. in section 4) and in the companion paper [1].

### 3.3. Probability distributions and simple sampling

To specify prior and data distributions in a Bayesian Problem we need implementations of probability distributions and CUQIpy provides the `Distribution` class for this.

While in general the data distribution should be chosen based on knowledge or assumptions about noise and/or the forward model, choosing a prior is somewhat more of a subjective task. Gaussian priors are often used, in part due to simplicity and because they often lead to computationally efficient sampling, but also because they provide a flexible framework for constructing expressive priors. One such example is GMRFs, which are often used to express smoothness in the prior [2], akin to Tikhonov regularization. On the other hand, Laplace and Cauchy distributions can be used as sparsity-inducing priors due to their heavy tails, while Laplace and Cauchy Markov random field (LMRF and CMRF) priors provide edge-preservation [2], similarly to TV regularization.

**Table 4.** Common distribution classes in CUQIpy. First block: simple distributions. Second block: composite distributions and Markov random fields. Third block: distributions used to specify Bayesian inverse problems. Fourth block: utility distribution classes providing access to benchmark and user-defined distributions.

| Distribution name | Description |
| --- | --- |
| `Beta` | Beta distribution |
| `Cauchy` | Cauchy distribution |
| `Gamma` | Gamma distribution |
| `Gaussian` | Multivariate Gaussian distribution |
| `InverseGamma` | Inverse Gamma distribution |
| `Laplace` | Laplace distribution |
| `Lognormal` | Log-normal distribution |
| `Uniform` | Uniform Distribution |
| `CMRF` | Cauchy Markov random field |
| `GMRF` | Gaussian Markov random field |
| `JointGaussianSqrtPrec` | Gaussian of multiple square-root precision matrices |
| `LMRF` | Laplace Markov random field |
| `JointDistribution` | Joint distribution of multiple parameters |
| `Posterior` | Posterior distribution |
| `DistributionGallery` | Collection of benchmark distributions |
| `UserDefinedDistribution` | Defined by functions for logarithm of probability density function `logpdf` and/or `sample` |

As of now, CUQIpy provides a selection of commonly used distributions, see table 4, including those mentioned above. Furthermore, the framework is general and it is easy to add new distributions. Users can also set up their own distributions through `UserDefinedDistribution` by providing a function to evaluate the logarithm of the PDF and/or a function to generate a single sample.

The rest of this section illustrates different choices of priors and how they perform when applied to a Bayesian inverse problem with a smooth signal and a piecewise constant signal, respectively.

To express an i.i.d. Gaussian prior on $\boldsymbol{x}$ with zero mean and standard deviation 0.1,

$$\boldsymbol{x} \sim \text{Gaussian}\left(\boldsymbol{0}, 0.1^2 \boldsymbol{I}_n\right). \tag{22}$$

The CUQIpy command would be

```
x = Gaussian(x_zero, 0.1**2, geometry=A.domain_geometry)
```

The last part equips the distribution with the domain geometry of $\boldsymbol{A}$; this is not required but will prove convenient.

Distributions can evaluate probability density function (`pdf`), cumulative density function (`cdf`), logarithm of pdf (`logpdf`) and gradient of the `logpdf` (`gradient`) for a particular element, e.g.

```
x.logpdf(x_zero)
```

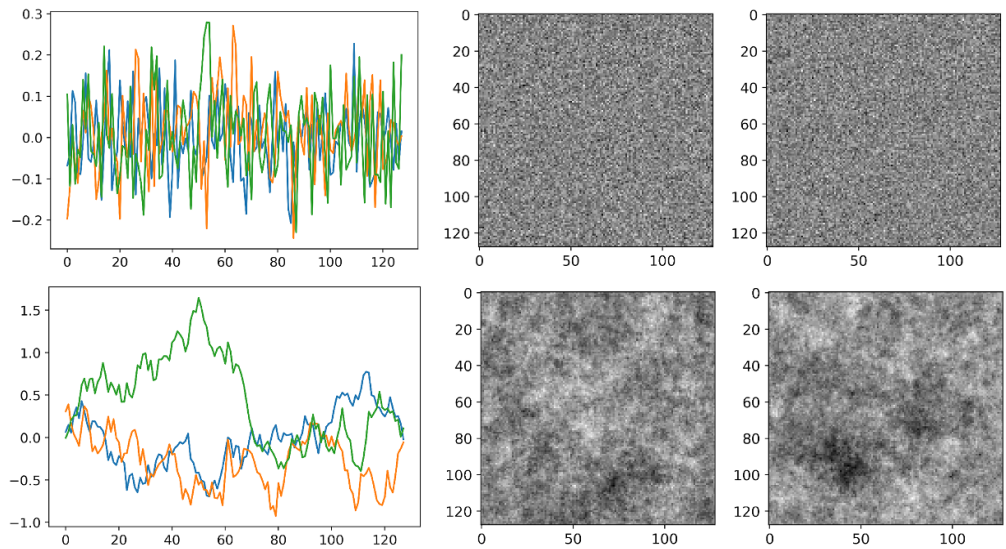

**Figure 4.** Top left: three samples from the Gaussian distribution for the 1D geometry. Top right: two samples for the 2D geometry. Bottom left: three samples from GMRF on 1D geometry. Bottom right: two samples for the 2D geometry.

Distributions can be sampled; for simple distributions this is done directly (i.e. without MCMC) using the `sample` method. Here we generate three samples and plot them:

```
samples_x = x.sample(3)
samples_x.plot()
```

Sampling returns a `Samples` object, which knows of the geometry from the distribution, so that the generated samples can automatically be displayed in the appropriate way, see figure 4 for 1D and 2D cases of the Gaussian distribution.

Another distribution example is a GMRF prior as discussed in section 2.3:

$$\boldsymbol{x} \sim \mathrm{GMRF}\left(\boldsymbol{0}, 50\right) \tag{23}$$

which can be specified both in 1D as

```
x = GMRF(x_zero, 50)
```

and in 2D by passing an `Image2D` geometry as in section 1.2; samples shown in figure 4.

As previously discussed, parameters such as a standard deviation may not be known and to be inferred along with the main parameter of interest. To set the stage for modelling of this in a later section we describe here conditional distributions in CUQIpy. For example in the Gaussian case, we can let $d$ be the unknown standard deviation

$$\boldsymbol{x} \sim \mathrm{Gaussian}\left(\boldsymbol{0}, d^2 \boldsymbol{I}_n\right) \tag{24}$$

by expressing this in CUQIpy using the syntax

```
x = Gaussian(x_zero, lambda d: d**2)
print(x)
```

Printing shows this is a conditional distribution, expressing that $x$ is conditional on $d$:

```
CUQI Gaussian. Conditioning variables ['d'].
```

Here a simple squared dependence on $d$ is specified via a lambda function; it is possible to express more involved dependencies in this way. We can now condition on specific values of $d$ to produce a fully specified distribution, from which sampling is possible, e.g.

```
x(d=0.1).sample(10)
x(d=0.3).sample(10)
```

produces 10 samples from $x|d$ for the two choices $d = 0.1$ and $d = 0.3$, respectively.

A common use case for conditional distributions is to express the *data distribution* for an inverse problem, in our example the conditional distribution $y|x$. Assuming Gaussian noise (with a known noise level) on the measured data, this distribution is given by

$$y \sim \text{Gaussian}\left(Ax, 0.01^2 I_m\right). \tag{25}$$

We can express this in CUQIpy by

```
y = Gaussian(A @ x, 0.01**2)
```

where the use of `@` emphasises the structural information that `A` is a linear operator. In general, lambda (anonymous) functions are used to define algebraic expressions on parameters, with some operations like `A @ x` being natively supported without lambda expressions. More native algebraic operations are planned in future versions of CUQIpy. In this scenario, the previously defined distribution x is employed to establish the dependency.

Calling `print(y)` shows that this is a conditional distribution on x:

```
CUQI Gaussian. Conditioning variables ['x'].
```

As we shall see in the following section this general syntax is used for specifying a Bayesian Problem to be solved for given observed data in a simple and intuitive way. Here, we show how one can simulate new data realizations by conditioning on the true signal, sampling the resulting distribution and plotting the resulting samples:

```
y(x=x_true).sample(5).plot()
```

These new example data realizations are shown in figure 5 top left.

### 3.4. Bayesian modelling: specifying inverse problems for UQ analysis

We now have all the tools we need to specify a complete Bayesian Problem in CUQIpy. We consider the 1D deblurring sinc case and choose a GMRF prior on the signal $x$ and assume Gaussian data distribution:

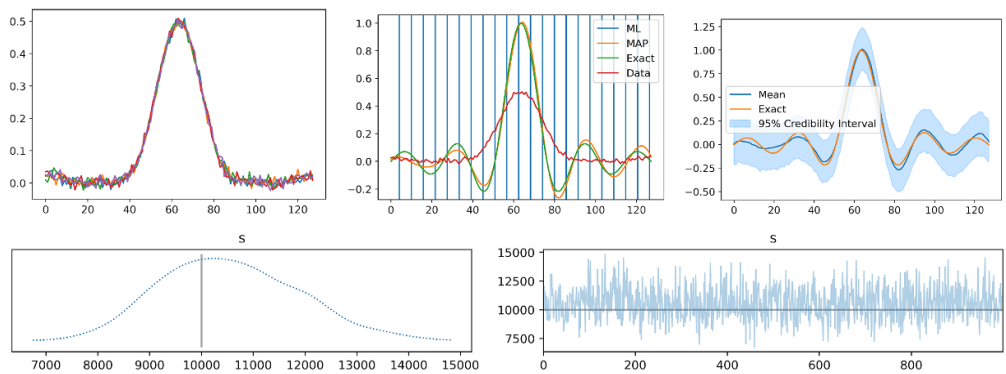

**Figure 5.** Results for Sinc example. Top left: samples of data distribution given $x^{\text{true}}$. Top center: ML, MAP, $x^{\text{true}}$ and $y^{\text{data}}$. Top right: posterior mean and 95% credibility interval compared to $x^{\text{true}}$. Bottom left: kernel density estimate from posterior samples of $s$. Bottom right: trace (chain) of posterior samples of $s$ compared to $s^{\text{true}}$.

$$x \sim \text{GMRF}(\mathbf{0}, 50), \qquad (26a)$$

$$y \sim \text{Gaussian}(Ax, 0.01^2 I_m) \qquad (26b)$$

which in CUQIpy is specified as

```
x = GMRF(x_zero, 50)
y = Gaussian(A @ x, 0.01**2)
```

The most automated high-level approach is to form a `BayesianProblem`:

```
BP = BayesianProblem(y, x)
print(BP)
```

Printing demonstrates how this represents a joint distribution $p(x,y)$:

```
BayesianProblem with target:
 JointDistribution(
    Equation:
        p(y,x) = p(y|x)p(x)
    Densities:
        y ~ CUQI Gaussian. Conditioning variables ['x'].
        x ~ CUQI GMRF.
)
```

Next, we specify the observed data $y^{\text{obs}}$ to produce a posterior distribution

```
BP.set_data(y=y_obs)
```

Maximum likelihood (ML) and maximum *a posteriori* (MAP) estimates are given:

```
x_ML = BP.ML()
x_MAP = BP.MAP()
x_ML.plot(), x_MAP.plot()
```

which are `CUQIarrays` containing the geometry so may be plotted directly, see figure 5.

We can also generate a desired number of posterior samples using MCMC by

```
samples = BP.sample_posterior(1000)
```

This will analyze the problem structure and try to determine a suitable sampler. In the present case of a linear forward model and Gaussian prior and data distribution (recall that GMRF is a Gaussian with special structure) the `LinearRTO` sampler is efficient and selected. This is possible for a range of problem types, with more being added. More details and how to manually specify a sampler will be described in section 3.6.

The MCMC samples are returned in a `Samples` object that as we have seen before offers various plotting options such as the posterior mean and credibility interval, with the option to compare with the exact:

```
samples.plot_ci(exact=x_true)
```

Finally, we mention the convenience 'UQ' method that carries out the above steps and generates selected solution plots:

```
BP.UQ()
```

### 3.5. Hierarchical modeling

Until now, we assumed a known noise level, i.e. the standard deviation on the noise in the data distribution. Often this is not known and we can then include this in the set of parameters to be estimated along with $x$, as mentioned in section 2.1. Such a hyperparameter is often modelled by a Gamma distribution on the precision (inverse variance), since it forms a conjugate pair with the Gaussian distribution, and this can be exploited for efficient sampling. The problem then takes the form

$$s \sim \text{Gamma}\left(1, 10^{-4}\right), \tag{27a}$$
$$x \sim \text{GMRF}(\mathbf{0}, 50), \tag{27b}$$
$$y \sim \text{Gaussian}\left(Ax, s^{-1} I_m\right) \tag{27c}$$

which we can specify in CUQIpy as (noting again the near-math syntax)

```
s = Gamma(1, 1e-4)
x = GMRF(x_zero, 50)
y = Gaussian(A @ x, prec=lambda s: s)
```

Note for illustration purposes we use the `prec` keyword to define the Gaussian with a precision rather than the default covariance used earlier. The same result can be achieved by using the syntax `cov = lambda s: 1/s`.

As before we set up a Bayesian Problem, this time also including *s*, and provide the observed data $y^{\text{obs}}$, leaving *x* and *s* to be estimated.

```
BP = BayesianProblem(y, x, s)
BP.set_data(y=y_obs)
```

We can then sample the posterior distribution; this will detect that multiple parameters are to be estimated and employ a Gibbs sampler, while making use of any conjugacy, here between the Gaussian and Gamma distributions, to employ efficient samplers for each of the conditional distributions.

```
samples = BP.sample_posterior(1000)
```

The sampling outputs the automatically determined sampling strategy selected for the posterior:

```
Using Gibbs sampler
Automatically determined sampling strategy:
    x: LinearRTO
    s: Conjugate
```

For a multi-parameter problem, the `Samples` object returned contains chains for each parameter that can be individually picked out and plotted, for example to plot the trace of *s* and the credibility interval for *x* compared to the ground truth `s_true = 1/0.01**2` (recall *s* was the precision) and `x_true` respectively, we call

```
samples["s"].plot_trace(exact=s_true)
samples["x"].plot_ci(exact=x_true)
```

The figure in figure 5 shows that the posterior for *s* centers around the known true value $s^{\text{true}} = 0.01^{-2}$ and the plot of individual samples has the 'fuzzy worm' appearance suggesting that samples are independent. It is also possible to compute the ESS, chain auto correlation and more via the `Samples` object.

We note that the 2D deconvolution example in the introduction (section 1) is set up and solved in precisely the same way. This highlights the powerful abstraction layer in CUQIpy that allows problems of different types, dimensions etc to be modelled and solved using almost identical code.

### 3.6. Manually choosing the sampling method

The class `BayesianProblem` is aimed at automating the sampling process allowing the user to focus on just the modelling. For full control, it is possible to set up and sample the problem manually. CUQIpy provides a range of samplers (table 5) that can be employed in case the automated sampling is not available or unsatisfactory. We demonstrate this for the square 1D deconvolution test problem with three different choices of prior.

**Table 5.** Samplers in CUQIpy grouped as they are presented in section 2.

| Sampler name | Description |
| --- | --- |
| MH | Metropolis–Hastings |
| pCN | preconditioned Crank–Nicolson |
| ULA | Unadjusted Langevin algorithm |
| MALA | Metropolis-Adjusted Langevin algorithm |
| NUTS | No U-Turn Sampler |
| LinearRTO | Linear Randomize-Then-Optimize |
| UGLA | Unadjusted Laplace Approximation |
| Gibbs | Gibbs sampler for joint distributions |
| CWMH | Component-Wise Metropolis–Hastings |
| Conjugate | Conjugate sampler |
| ConjugateApprox | Approximate conjugate sampler |

For the smooth sinc signal, the GMRF produced accurate MAP and posterior mean solutions, but for non-smooth square signal we expect edge-preserving priors to be more suitable. We compare GMRF with LMRF and CMRF by letting

```
x = GMRF(x_zero, prec=50)             # or
x = LMRF(x_zero, scale=0.01)          # or
x = CMRF(x_zero, scale=0.01)
```

in combination with

```
y = Gaussian(A @ x, 0.01**2)
```

As was seen from the output of printing a Bayesian Problem earlier, it contains a target *joint distribution* over all the parameters of interest. To sample manually, instead of setting up the Bayesian Problem we set up this joint distribution over the *x* and *y* parameters, and—being a distribution—we condition it on the observed data $y^{\text{obs}}$ to obtain the posterior distribution:

```
joint = JointDistribution(y, x)
posterior = joint(y=y_obs)
```

Each of the three priors have a certain structure that can be exploited for efficient sampling by means of a choosing a suitable sampler

```
mySampler = LinearRTO(posterior)      # For GMRF
mySampler = UGLA(posterior)           # For LMRF
mySampler = NUTS(posterior)           # For CMRF
```

As already mentioned the `LinearRTO` sampler is a good choice for the GMRF case. For the LMRF case we use a dedicated sampler that uses a Gaussian approximation of the prior for efficiency known as the unadjusted Laplace approximation [58] and denoted as `UGLA` in CUQIpy. In case of CMRF, we exploit that the density function is differentiable and employ the gradient-based sampler NUTS. We note that these choices of samplers exactly match what is done by BayesianProblem for GMRF, LMRF, and CMRF priors, respectively.

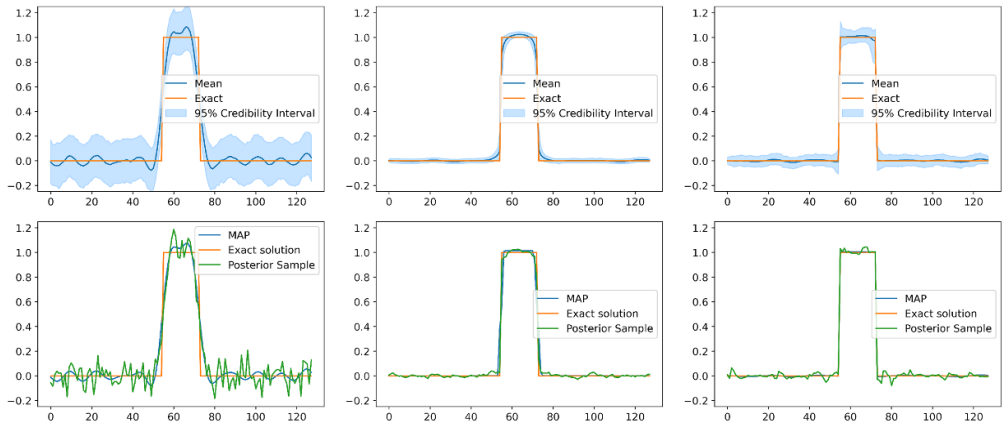

**Figure 6.** UQ results for the square signal. Left: GMRF prior for this test problem. Center: LMRF prior. Right: CMRF prior. Top row: posterior mean and 95% credibility interval compared to exact solution. Bottom row: MAP estimate and a posterior sample compared to exact solution.

Having specified the posterior and a sampler for each choice of prior, we proceed to generate a desired number of samples (1000), after a burn-in phase of 200 samples:

```
samples = mySampler.sample(1000, 200)
```

after which posterior information can be plotted:

```
samples.plot_ci(exact=x_true)
```

The MAP and ML estimates can also be manually computed via the solver module of CUQIpy, but this is omitted for brevity.

The results of the UQ analysis are shown in figure 6. Clearly the GMRF prior (which imposes smoothness) is not suitable for the piecewise constant signal: it produces a smooth posterior mean solution with large variability. We also plot a single sample (the 500th) which confirms large variability. The LMRF prior and especially the CMRF prior produce solutions that much better resembles a piecewise constant signal, in particular the MAP solution, and much less variability.

The samplers also return diagnostic information into the samples object, as demonstrated in section 4. A number of plotting methods are available in the samples objects; some are already demonstrated, and some will be illustrated in the following case studies.

Having demonstrated some of the key concepts of CUQIpy on a simple test problem, we now present three more complex case studies designed to highlight the capabilities and flexibility to model, solve and conduct UQ analysis for inverse problems. Further examples involving PDE-based inverse problems are given in the companion paper [1].

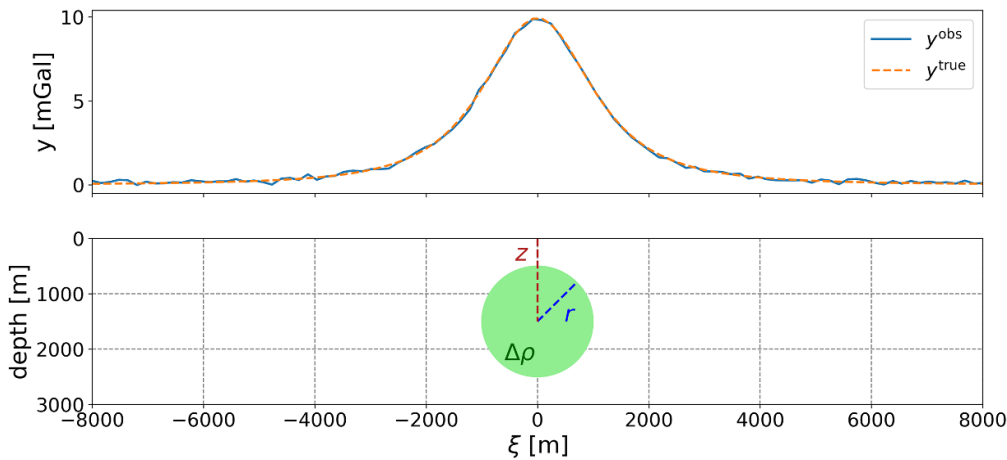

**Figure 7.** Top: simulated measured gravity anomaly with noise and true gravity anomaly. Bottom: buried spherical body causing the gravity anomaly at the surface.

## 4. Case study: user-specified nonlinear model with a gravity anomaly

In this case study we demonstrate how CUQIpy can be applied to solve a non-linear inverse problem. We focus on an example from geophysics, where a subsurface density contrast gives rise to a gravity anomaly. This problem is not included in CUQIpy's suite of test problems, and hence it illustrates how a user-specified model can be handled in CUQIpy. We also demonstrate how different samplers can easily be applied to the same problem, and we will see that the NUTS sampler is much better suited for sampling correlated parameters than a simple MH sampler. Below, we specify the model and the Bayesian Problem using tools imported from CUQIpy:

```
from cuqi.model import Model
from cuqi.distribution import Gaussian, JointDistribution
from cuqi.sampler import MH, NUTS
from cuqi.geometry import Continuous1D, Discrete
```

### 4.1. User-specified forward model of gravity anomaly measurements

A subsurface body with a density contrast to its homogeneous surroundings causes a gravity anomaly. The gravity anomaly field depends on the body's depth, shape, size, and density contrast. We construct an example where a spherical body causes the gravity anomaly shown in figure 7. The sphere has radius $r^{\text{true}} = 1000\,\text{m}$, the density contrast is $\Delta\rho^{\text{true}} = 800\,\text{kg}\,\text{m}^{-3}$, and the depth from surface to center of the sphere $z^{\text{true}} = 1500\,\text{m}$. This could be interpreted as a crude model of iron ore (a mineral substance with high iron content) buried in sedimentary rock.

In the inverse problem we are interested in inferring the spherical body's radius, density contrast and depth given measurements of the gravity anomaly at the surface. We measure the vertical component of the gravity anomaly, caused by the spherical body, at the surface along a $\xi$-axis whose origin intersects the point above the center of the sphere, see figure 7. The measured signal at the position $\xi$ is modelled as arising from the parameters of interest $z, \Delta\rho, r$ through the non-linear function

$$f(z, \Delta\rho, r, \xi) = \frac{4\pi}{3} G \left(\frac{\Delta\rho r^3}{z^2}\right) \left(\frac{1}{1+(\xi/z)^2}\right)^{3/2}, \tag{28}$$

where $G$ is the gravitational constant. We measure the signal at $m$ points $\xi_1, \ldots, \xi_m$ to obtain the data vector $\boldsymbol{y} \in \mathbb{R}^m$ with elements given by

$$y_i = f(z, \Delta\rho, r, \xi_i) = A_i(\boldsymbol{x}), \qquad \boldsymbol{x} = [z, \Delta\rho, r]^\mathsf{T}, \tag{29}$$

and hence the complete nonlinear forward model $\boldsymbol{A} : \mathbb{R}^3 \to \mathbb{R}^m$ is expressed as

$$\boldsymbol{y} = \boldsymbol{A}(\boldsymbol{x}) = \begin{bmatrix} A_1(\boldsymbol{x}) \\ \vdots \\ A_m(\boldsymbol{x}) \end{bmatrix}. \tag{30}$$

The forward model is implemented in `forward_gravity()`, cf appendix. To use gradient-based samplers we need derivative information of the forward model, i.e. the Jacobian matrix $\boldsymbol{J} \in \mathbb{R}^{m \times 3}$ of the forward operator, with elements defined by

$$J_{ij} = \frac{\partial A_i(\boldsymbol{x})}{\partial x_j}. \tag{31}$$

Computation of the Jacobian is implemented in `jac_gravity()`, cf appendix.

The user-specified CUQIpy model must contain a python function that evaluates the forward model as well as information about the problem's range and domain geometries. The measurements $\boldsymbol{y}$ are given by the continuous function $f$ in (28) while the inferred parameters in $\boldsymbol{x}$ consist of three different physical quantities, making the CUQIpy geometries `Continuous1D()` and `Discrete()` appropriate for the range and domain, respectively. Derivative information is added via a python function that evaluates the Jacobian matrix. The CUQIpy model is specified as follows:

```
A = Model(forward=forward_gravity,
          jacobian=jac_gravity,
          range_geometry=Continuous1D(m),
          domain_geometry=Discrete(["z", "rho", "r"]))
```

### 4.2. The Bayesian problem

We specify a simple Gaussian prior with means that are different from the true parameter values and with large standard deviations, such that the prior is not very informative. Furthermore, we assume a Gaussian data distribution with standard deviation $10^{-6}$ m s$^{-2}$. This leads to the Bayesian Problem:

$$\boldsymbol{x} \sim \text{Gaussian}\left([1550, 850, 950]^\mathsf{T}, \text{diag}\left(500^2, 300^2, 300^2\right)\right), \tag{32a}$$

$$\boldsymbol{y} \sim \text{Gaussian}\left(\boldsymbol{A}(\boldsymbol{x}), \left(10^{-6}\right)^2 \boldsymbol{I}_m\right), \tag{32b}$$

where 'diag' denotes a diagonal matrix or in **CUQIpy**:

```
x = Gaussian(np.array([1550, 850, 950]),
             sqrtcov=np.array([500, 300, 300]))
y = Gaussian(A(x), sqrtcov=1e-6)
```

To solve the inverse problem, we simulate a gravity anomaly signal $\boldsymbol{y}^{\text{obs}}$ by obtaining one realization from the data distribution (32*b*) at $m = 100$ equidistantly distributed points between $\xi = -8000$ m and $\xi = 8000$ m, and with $\boldsymbol{x}^{\text{true}} = [z^{\text{true}}, \Delta\rho^{\text{true}}, r^{\text{true}}]^{\mathsf{T}}$. In **CUQIpy** this is conveniently done using the data distribution's `sample()` method:

```
y_obs = y(x=x_true).sample()
```

The true and the noisy signals are plotted in figure 7.

### 4.3. Posterior sampling and analysis

In **CUQIpy**, we define the posterior distribution by forming the joint distribution of the prior and data distributions and conditioning on the observed data:

```
posterior = JointDistribution(x, y)(y=y_obs)
```

Examining (28) we see that the signal is determined from $z$ and the product $\Delta\rho r^3$, which means that $z$ can be inferred from the data while $\Delta\rho$ and $r$ can not be resolved individually. Therefore, we expect that the $\Delta\rho$ and $r$ samples will be extremely correlated. We compare the performance of the MH sampler and the NUTS sampler. We use the same starting point $\boldsymbol{x}^{\text{init}} = [1000, 2000, 1000]^{\mathsf{T}}$ and run both samplers for approximately 10 min.

We run the MH sampler with adaptive choice of step size by:

```
MHsampler = MH(posterior, x0=x_init, scale=100)
samplesMH = MHsampler.sample_adapt(1000000, 100000)
```

and the NUTS sampler with:

```
NUTSsampler = NUTS(posterior, x0=x_init)
samplesNUTS = NUTSsampler.sample_adapt(8000, 100)
```

We analyse the samples by inspecting correlation coefficients between the inferred variables and trace plots of the samples using **CUQIpy**'s plotting functionality, e.g. `samplesMH.plot_trace()`. As expected, table 6 reveals strong correlation between the radius $r$ and density contrast $\Delta\rho$ samples, whereas the depth $z$ samples are not correlated with the other parameters. This is the case for both sampling methods.

The trace plots in figure 8 show that for both sampling methods the depth parameter chains are well-mixed and that the posterior 1D marginals converge to the same distribution. More interestingly, the trace plots show that the MH chains of density contrast and radius are very poorly mixed, and the posterior 1D marginals have not converged to distributions of a simple form. The poor mixing is also confirmed by the ESS computed with **CUQIpy**'s `samplesMH.compute_ess()`. On the other hand, the NUTS sampler, that utilizes gradient information, immensely improves the mixing of the chains [37]. This is indicated by the chain plots, the smoother 1D marginals, and that the much larger ESS for NUTS, even though the number of posterior samples is much smaller. This demonstrates that gradient-based samplers like NUTS are better suited for correlated parameters.

**Table 6.** Correlation coefficients between sampled parameters using (left) the Metropolis–Hastings sampler and (right) the NUTS sampler.

|  | $z$ | $\Delta\rho$ |  |  | $z$ | $\Delta\rho$ |
|---|---|---|---|---|---|---|
| $\Delta\rho$ | 0.004 |  |  | $\Delta\rho$ | $-0.01$ |  |
| $r$ | 0.03 | $-0.98$ |  | $r$ | 0.03 | $-0.95$ |

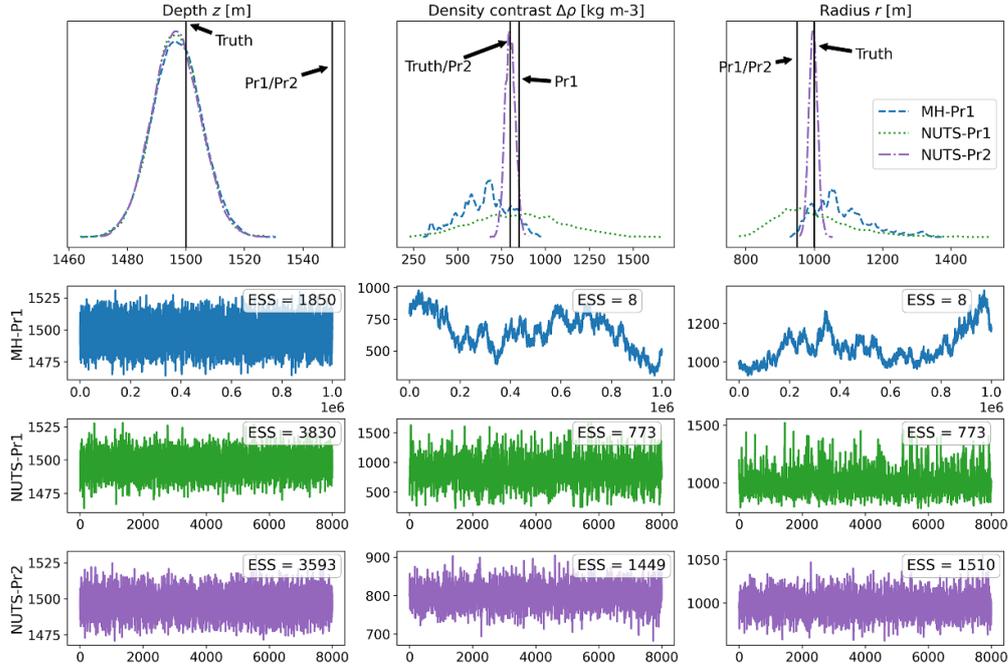

**Figure 8.** Trace plots of posterior samples using the uninformative prior (Pr1) in combination with the Metropolis–Hastings and NUTS samplers, and the strong $\Delta\rho$ prior (Pr2) in combination with the NUTS sampler. The top row show the 1D marginals of the posterior distributions with the prior and true parameter values marked with black vertical lines. Rows 2–4 show the chains of the samples and respective ESS that are computed for each chain using `samples.compute_ess()`.

Finally, we see that the depth posterior marginal seems to contract close to the true value with fairly high certainty, and is not very biased towards the prior mean. On the other hand, the posterior marginals of the density contrast and the radius contract around their respective prior means, but with high uncertainty. This reflects the fact that the data informs the depth posterior well, while the density contrast and radius are more controlled by the prior, as the data can not resolve these parameters individually.

To improve the posterior estimates we formulate a Gaussian prior that is more informative about the density contrast, but less informative about the radius or depth. This simulates a scenario where we seek iron ore in sedimentary rock with unknown size and depth. Hence, the depth and radius priors remain unchanged but the density contrast prior now has the true density contrast as mean, and a small standard deviation. The new, more informative prior is:

$$x \sim \text{Gaussian}\left([1550, 850, 950]^\top, \text{diag}\left(500^2, 30^2, 300^2\right)\right). \quad (33)$$

Since the density contrast and radius remain correlated, we use NUTS to sample the posterior as before. Using the new $\Delta\rho$ prior, we see in figure 8 that both the $\Delta\rho$ posterior and the $r$ posterior contract around the ground truth with high certainty.

## 5. Case study: x-ray CT with the CUQIpy-CIL plugin

In this case study we demonstrate how CUQIpy can be applied to a larger real-data 2D imaging problem in x-ray CT through the CUQIpy-CIL plugin that provides a simple interface to configure a model for CT. The only component needed to allow us to handle a CT problem in CUQIpy is the forward model which we load from the plugin, along with general CUQIpy tools:

```
from cuqipy_cil.model import FanBeam2DModel
from cuqi.array import CUQIarray
from cuqi.geometry import Image2D
from cuqi.distribution import Gaussian, Gamma, LMRF, JointDistribution
from cuqi.sampler import UGLA, Conjugate, ConjugateApprox, Gibbs
```

We assume the data and prior distributions are known, but with unknown precision parameters, and we show how this is modelled with a simple hierarchical model. Finally, we sample from the posterior and display the results via CUQIpy plotting functionality.

We use dataset B from the Helsinki Tomography Challenge (HTC) 2022 [36] to demonstrate UQ for CT. This is a 2D dataset and it is measured using fan-beam geometry. The data comes with a set of parameters characterizing the measurement setup and we use them to construct the model in section 5.1. To make the data compatible with CUQIpy we load it into a CUQIarray and equip it with geometry.

```
y_obs = CUQIarray(sinogram, geometry=Image2D(im_shape=(360, 560)))
```

This also allows easy data plotting by `y_obs.plot()`, see figure 9(a).

### 5.1. Model

In CT, an image of an object's interior is found from observed attenuation of x-rays passing through the object [23]; this is often formulated as a linear inverse problem $y = Ax$. Here, $x$ represents an image of the object's spatially varying attenuation coefficient, $y$ represents the measured attenuation (sometimes called the sinogram), and each matrix element $A_{ij}$ represents the $i$th x-ray's intersection with pixel $j$.

We use the CUQIpy-CIL plugin to specify the CT model. It is based on the Python package CIL [28, 39] from which the forward and adjoint operations (projection and backprojection) are wrapped as a CUQIpy model. Furthermore CUQIpy-CIL sets up the model with Image2D domain and range geometries. The CUQIpy-CIL model is straightforward to specify using parameters from the HTC dataset:

```
A = FanBeam2DModel(det_count=560,
                   det_spacing=0.2,
                   angles=-np.linspace(0, 2*np.pi, 360),
                   source_object_dist=410.66,
                   object_detector_dist=143.08,
                   domain=(83.09, 83.09),
                   im_size=(500,500))
```

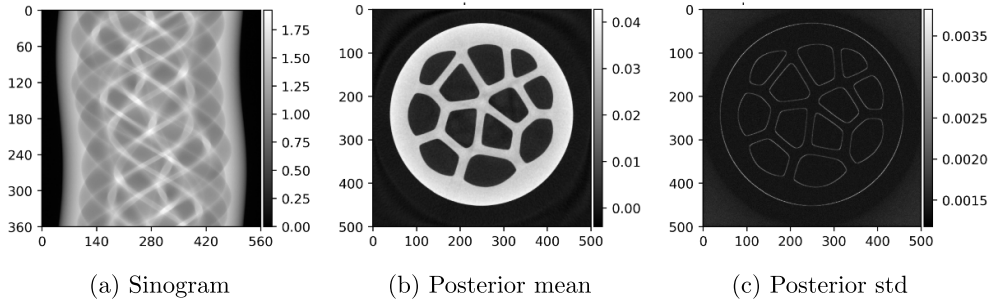

(a) Sinogram　　　　(b) Posterior mean　　　　(c) Posterior std

**Figure 9.** (a) Sinogram containing the noisy data $y^{\text{obs}}$. (b)–(c) Visualization of some basic statistics of the $x$ posterior samples.

### 5.2. The Bayesian problem

To construct our Bayesian Problem we must specify the data and prior distributions. We want an edge-preserving prior because we expect a piecewise constant image. Therefore, we apply a 2D LMRF prior to the unknown image $x$, which is similar to total-variation regularization [22, section 8.6]. Since we do not know the prior precision, we let this be a random variable $d$ that follows a Gamma distribution, and we include this in the prior. The data distribution is formulated assuming the data noise is Gaussian, which is a good approximation for high-intensity x-rays. However, we do not know the noise level of the dataset. Therefore we model the noise precision as a random variable $s$ and let it follow a Gamma distribution.

With these prior and data distributions, we define the Bayesian Problem:

$$d \sim \text{Gamma}\left(1, 10^{-4}\right), \tag{34a}$$

$$s \sim \text{Gamma}\left(1, 10^{-4}\right), \tag{34b}$$

$$x \sim \text{LMRF}\left(\mathbf{0}, d^{-1}\right), \tag{34c}$$

$$y \sim \text{Gaussian}\left(Ax, s^{-1}\right). \tag{34d}$$

In CUQIpy we express this by:

```
d = Gamma(1, 1e-4)
s = Gamma(1, 1e-4)
x = LMRF(0, lambda d: 1/d, geometry=A.domain_geometry)
y = Gaussian(A @ x, lambda s: 1/s)
```

*5.3. Posterior sampling and analysis*

We set up the posterior as a joint distribution and condition on the observed data:

```
posterior = JointDistribution(d, s, x, y)(y=y_obs)
```

This posterior has a hierarchical structure so we use a Gibbs sampler. We utilize that the hyperpriors for the precision parameters *l* and *d* are conjugate to the data and prior distributions respectively. In CUQIpy we construct the appropriate Gibbs sampler, and we obtain 500 samples from the posterior after a burn-in of 100 samples:

```
sampling_strategy = {'d': ConjugateApprox,
                     's': Conjugate,
                     'x': UGLA}
samples = Gibbs(posterior, sampling_strategy).sample(500, 100)
```

These match the automatic sampler selection done by `BayesianProblem`.

We analyze and visualize posterior samples and basic statistics for the CT reconstruction $x$ and the hyperparameters *s* and *d*:

```
samples["x"].plot_mean() ,    samples["x"].plot_std()
samples["s"].plot_trace(),    samples["d"].plot_trace()
```

The posterior mean for $x$ is a good reconstruction of the CT data, and its standard deviation image highlights most uncertainty near the edges (figure 9). The trace plots for *s* and *d* in figure 10 suggest that chains have converged and we see that the distributions of the samples *s* and *d* take their maxima around 9250 and 460. These parameters estimate the inferred noise and regularization levels and how certain the estimates are.

## 6. Case study: AD benchmark by CUQIpy-PyTorch

In this case study, we demonstrate the capabilities of the CUQIpy-PyTorch plugin, which enables the use of PyTorch's AD framework within CUQIpy, allowing efficient gradient-based sampling of an arbitrary Bayesian Problem using the NUTS sampler. To showcase this, we use a benchmark problem that is commonly used to illustrate Bayesian sampling with NUTS, known as the Eight Schools problem.

The Eight Schools problem is a classic example in Bayesian statistics, first introduced in [49] and later adapted as a benchmark problem in STAN [10], Pyro [5], Edward [57], and

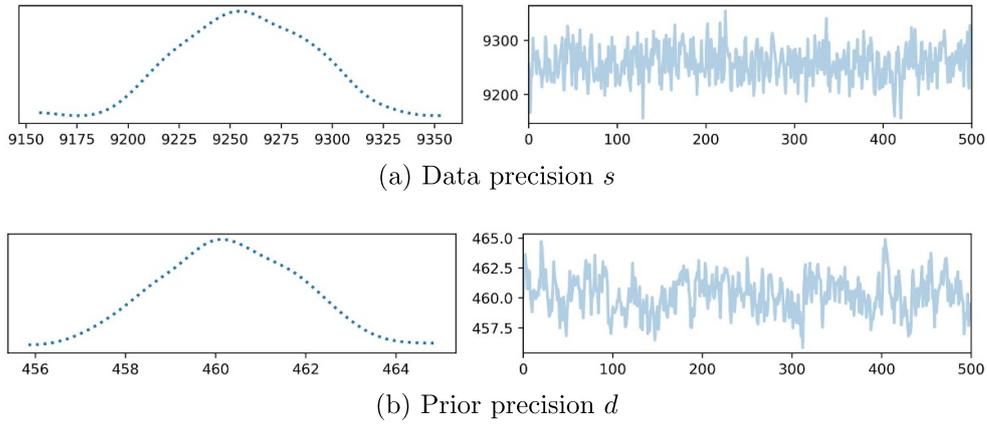

(a) Data precision $s$

(b) Prior precision $d$

**Figure 10.** Trace plots of posterior samples of the hyperparameters. Left column: posterior distributions. Right column: sample chains.

Tensorflow probability [43]. The problem involves estimating the effects of coaching on exam scores at eight different schools. Specifically, we want to estimate the treatment effect, $x_i$, at each school $i$, which represents the difference in exam scores between students who received coaching and those who did not. However, because the number of students in each school is small, the estimates of the treatment effects are noisy and may be biased. We therefore use a Bayesian hierarchical model to estimate the treatment effects, which accounts for both within-school variation and between-school variation. For more details see [16, 49] and the above-mentioned software packages.

Treatment effects and standard deviations are estimated in a study and assumed to be observed and known. We import the packages and define the observations by:

```
import torch as xp
from cuqi.distribution import JointDistribution
from cuqipy_pytorch.distribution import Gaussian, LogGaussian
from cuqipy_pytorch.sampler import NUTS

# Observed treatment effect y_obs and standard deviations s_obs.
y_obs = xp.tensor([28,  8, -3,  7, -1,  1, 18, 12], dtype=xp.float32)
s_obs = xp.tensor([15, 10, 16, 11,  9, 11, 10, 18], dtype=xp.float32)
```

The centered hierarchical Bayesian Problem for Eight Schools is given by:

$$u \sim \text{Gaussian}\left(0, 10^2\right), \qquad (35a)$$
$$t \sim \text{Lognormal}(5, 1), \qquad (35b)$$
$$\boldsymbol{x}' \sim \text{Gaussian}\left(\boldsymbol{0}, \boldsymbol{I}_8\right), \qquad (35c)$$
$$x_i = u + t x'_i \quad i = 1, \ldots, 8, \qquad (35d)$$
$$\boldsymbol{y} \sim \text{Gaussian}\left(\boldsymbol{x}, \text{diag}\left(\boldsymbol{s}^{\text{obs}}\right)^2\right), \qquad (35e)$$

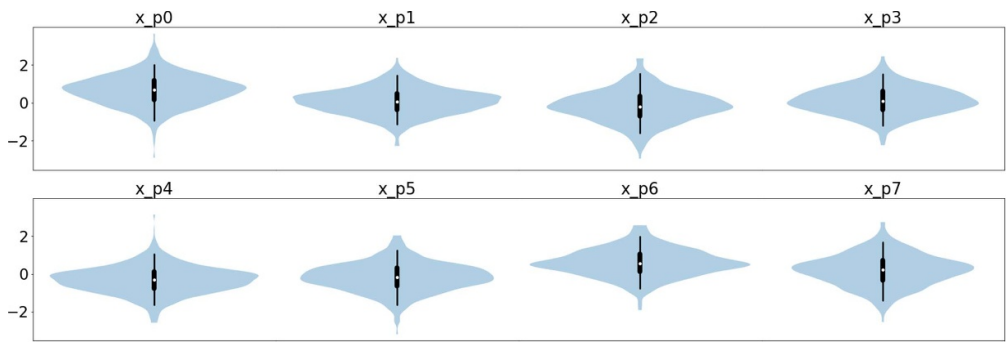

**Figure 11.** Violin plots showing combined kernel density estimates and box plots of the estimated standardized effects $x'$ obtained via NUTS sampling of the joint posterior. The results match those obtained with other software packages for the same problem and show a difference in the effectiveness of the coaching for the eight schools studied.

where $x = [x_1, \ldots, x_8]^\mathsf{T}$ represents the treatment effects at each school, $x' = [x'_1, \ldots, x'_8]^\mathsf{T}$ represents the standardized effects, $u$ represents the average treatment effect across all schools, $t$ represents the standard deviation of the treatment effects, and $y = [y_1, \ldots, y_8]^\mathsf{T}$ represents the exam score differences. Note that due to parameter relations there are no relations between distributions that can be exploited, e.g. for conjugacy. The Bayesian Problem can be implemented with CUQIpy-PyTorch as follows

```
u   = Gaussian(0, 10**2)
t   = Lognormal(5, 1)
x_p = Gaussian(xp.zeros(8), 1)
x   = lambda u, d, x_p: u+t*x_p
y   = Gaussian(x, s_obs**2)
```

Note here the use of a lambda function to define a function representing the random variable $x$. Since $x$ is entirely specified by the other random variables, the posterior of interest is $p(u,t,x'|y=y^{\text{obs}})$. We sample the problem by the optimized NUTS from Pyro [5], which is embedded in CUQIpy-PyTorch. We draw 1000 posterior samples and 500 warm-up samples with default target acceptance rate of 0.8:

```
joint = JointDistribution(u, t, x_p, y)   # Define joint distribution
posterior = joint(y=y_obs)                 # Create posterior
sampler = NUTS(posterior)                  # Define sampling strategy
samples = sampler.sample(N=1000, Nb=500)   # Sample from posterior
```

Because of the AD capabilities of CUQIpy-PyTorch the posterior is jointly sampled with NUTS. The analysis of the Eight Schools problem can be visualized, e.g. with a so-called violin plot, which is also available in CUQIpy, see figure 11. The sample means of the average treatment effect $u$ is found to be 5.6 with standard deviation $t$ of 12.7.

```
# Plot posterior samples
samples["x_p"].plot_violin();
print(samples["u"].mean()) # Average effect
print(samples["t"].mean()) # Average standard deviation
```

These results can be used to conclude on the effect of coaching for SAT scores [16]. The results obtained are in line with those obtained by the above-mentioned software packages.

## 7. Conclusion and future work

In this paper, we have presented CUQIpy, a Python package designed to serve as a versatile framework for UQ in inverse problems. This flexible framework enables users to focus on modeling aspects while still providing experts with access to the underlying 'machine room'.

We demonstrated the capabilities of CUQIpy by showcasing Bayesian inverse problems using both simple and advanced distributions, along with sampling strategies ranging from classical samplers to advanced techniques incorporating gradient information, forward model linearity, or distribution approximations.

As CUQIpy remains under active development, we anticipate further enhancement and expansion of its features beyond the current version 1.0.0. Future work for CUQIpy includes:

- Extending CUQIpy with new distributions and samplers.
- Extending support for automated sampler selection by `BayesianProblem` to accommodate a broader range of problems, including more hierarchical models.
- Enhancing the exploitation of structure through a wider range of algebraic operations on random variables and forward models, replacing lambda functions.
- Improving interoperability with other libraries and array types, inspired by NEP 47 (https://numpy.org/neps/nep-0047-array-api-standard.html).
- Incorporating support for UQ techniques beyond MCMC-based sampling, such as Variational Inference [6], to enable efficient sampling of posterior distribution approximations.
- Incorporate support for learned models and priors such as Plug & Play priors [31].

The flexibility and adaptability of CUQIpy encompass a wide array of inverse problems illustrating its potential as a powerful tool for UQ in the field. It is our hope that CUQIpy and similar computational tools will foster a more wide-spread adoption of UQ methods for large-scale inverse problems.

## Data availability statement

CUQIpy and plugins are available from https://cuqi-dtu.github.io/CUQIpy. The code and data to reproduce the results and figures of the present paper are available from https://github.com/CUQI-DTU/paper-CUQIpy-1-Core.

The data that support the findings of this study are openly available at the following URL/DOI: https://zenodo.org/doi/10.5281/zenodo.10512533.

## Acknowledgments

This work was supported by The Villum Foundation (Grant No. 25893). J S J would like to thank the Isaac Newton Institute for Mathematical Sciences for support and hospitality during the programme 'Rich and Nonlinear Tomography—a multidisciplinary approach' when work on this paper was undertaken. This work was supported by EPSRC Grant Number EP/R014604/1. This work was partially supported by a grant from the Simons Foundation (J S J). F U has been supported in part by Academy of Finland (Project Number 353095). The authors are grateful to members of the CUQI project for valuable input that helped shape the design of CUQIpy.

## Appendix. Python implementations for case study 3: gravity anomaly

This Python function implements the gravity model (28) evaluated in a point `wrt`:

```python
# Measurement grid
xi = np.linspace(-8000, 8000, 100)
# Gravitational constant N m^2 / kg^2
G = 6.6743e-11
# Forward function
def forward_gravity(wrt):
    z = wrt[0]
    rho = wrt[1]
    r = wrt[2]
    y = 4/3*np.pi*G*(rho*r**3/z**2)*(1/(1+(xi/z)**2))**(3/2)
    return y
```

This Python function evaluates the Jacobian of the gravity model in a point `wrt`:

```python
# Jacobian of forward function
def jac_gravity(wrt):
    z = wrt[0]
    rho = wrt[1]
    r = wrt[2]
    dAdz = -4/3*np.pi*G*rho*r**3*z* (2*z**2-xi**2)
            * 1/(z**2/(xi**2+z**2))**(1/2) * 1/(xi**2+z**2)**3
    dAdrho = 4/3*np.pi*G*(r**3/z**2)*(1/(1+(xi/z)**2))**(3/2)
    dAdr = 4*np.pi*G*(rho*r**2/z**2)*(1/(1+(xi/z)**2))**(3/2)
    J = np.vstack([dAdz, dAdrho, dAdr]).T
    return J
```

## ORCID iDs

Nicolai A B Riis  https://orcid.org/0000-0002-6883-9078
Amal M A Alghamdi  https://orcid.org/0000-0003-0145-5296

Felipe Uribe  https://orcid.org/0000-0002-1010-8184
Silja L Christensen  https://orcid.org/0000-0003-3995-3055
Babak M Afkham  https://orcid.org/0000-0003-3203-8874
Per Christian Hansen  https://orcid.org/0000-0002-7333-7216
Jakob S Jørgensen  https://orcid.org/0000-0001-9114-754X

## References

[1] Alghamdi A M A, Riis N A B, Afkham B M, Uribe F, Christensen S L, Hansen P C and Jørgensen J S 2024 CUQIpy: II. Computational uncertainty quantification for PDE-based inverse problems in Python *Inverse Problems* **40** 045010

[2] Bardsley J M 2019 *Computational Uncertainty Quantification for Inverse Problems* (SIAM)

[3] Bardsley J M, Solonen A, Haario H and Laine M 2014 Randomize-then-optimize: a method for sampling from posterior distributions in nonlinear inverse problems *SIAM J. Sci. Comput.* **36** A1895–910

[4] Biegler L, Biros G, Ghattas O, Heinkenschloss M, Keyes D, Mallick B, Marzouk Y, Tenorio L, van Bloemen Waanders B and Willcox K (eds) 2010 *Large-Scale Inverse Problems and Quantification of Uncertainty* (Wiley)

[5] Bingham E, Chen J P, Jankowiak M, Obermeyer F, Pradhan N, Karaletsos T, Singh R, Szerlip P, Horsfall P and Goodman N D 2019 Pyro: deep universal probabilistic programming *J. Mach. Learn. Res.* **20** 973–8 (available at: http://jmlr.org/papers/v20/18-403.html)

[6] Blei D M, Kucukelbir A and McAuliffe J D 2017 Variational inference: a review for statisticians *J. Am. Stat. Assoc.* **112** 859–77

[7] Brooks S, Gelman A, Jones G L and Meng X-L (eds) 2011 *Handbook of Markov Chain Monte Carlo* (CRC Press)

[8] Buzug T M 2008 *Computed Tomography* (Springer)

[9] Calvetti D and Somersalo E 2007 *Introduction to Bayesian Scientific Computing* (Springer)

[10] Carpenter B, Gelman A, Hoffman M D, Lee D, Goodrich B, Betancourt M, Brubaker M, Guo J, Li P and Riddell A 2017 Stan: a probabilistic programming language *J. Stat. Softw.* **76** 1–32

[11] Chung J, Knepper S and Nagy J G 2011 Large-scale inverse problems in imaging *Handbook of Mathematical Methods in Imaging* ed O Scherzer (Springer)

[12] Cotter S L, Roberts G O, Stuart A M and White D 2013 MCMC methods for functions: modifying old algorithms to make them faster *Stat. Sci.* **28** 424–46

[13] Dashti M and Stuart A M 2017 The Bayesian approach to inverse problems *Handbook of Uncertainty Quantification* ed R Ghanem, D Higdon and H Owhadi (Springer) ch 10, pp 311–428

[14] Gamerman D and Lopes H F 2006 *Markov Chain Monte Carlo: Stochastic Simulation for Bayesian Inference* 2nd edn (Chapman & Hall/CRC)

[15] Ge H, Xu K and Ghahramani Z 2018 Turing: a language for flexible probabilistic inference *Int. Conf. on Artificial Intelligence and Statistics (AISTATS 2018) (Playa Blanca, Lanzarote, Canary Islands, Spain, 9–11 April 2018)* pp 1682–90

[16] Gelman A, Carlin J B, Stern H S, Dunson D B, Vehtari A and Rubin D B 2013 *Bayesian Data Analysis* (CRC Press)

[17] Gelman A, Vehtari A, Simpson D, Margossian C C, Carpenter B, Yao Y, Kennedy L, Gabry J, Bürkner P-C and Modrák M 2020 Bayesian workflow (arXiv:2011.01808)

[18] Geman S and Geman D 1984 Stochastic relaxation, Gibbs distributions and the Bayesian restoration of images *IEEE Trans. Pattern Anal. Mach. Intell.* **PAMI-6** 721–41

[19] Ghattas O and Willcox K 2021 Learning physics-based models from data: perspectives from inverse problems and model reduction *Acta Numer.* **30** 445–554

[20] Groetsch C W 1993 *Inverse Problems in the Mathematical Sciences* (Vieweg)

[21] Hairer M, Stuart A M, Voss J and Wiberg P 2005 Analysis of SPDEs arising in path sampling part I: the Gaussian case *Commun. Math. Sci.* **3** 587–603

[22] Hansen P C 2010 *Discrete Inverse Problems: Insight and Algorithms* (SIAM)

[23] Hansen P C, Jørgensen J S and Lionheart W R B (eds) 2021 *Computed Tomography: Algorithms, Insight and Just Enough Theory* (SIAM)

[24] Hansen T M, Cordua K S, Looms M C and Mosegaard K 2013 SIPPI: a Matlab toolbox for sampling the solution to inverse problems with complex prior information: part 1—methodology *Comput. Geosci.* **52** 470–80
[25] Harris C R *et al* 2020 Array programming with NumPy *Nature* **585** 357–62
[26] Hastings W K 1970 Monte Carlo sampling methods using Markov chains and their applications *Biometrika* **57** 97–109
[27] Hoffman M D and Gelman A 2014 The No-U-Turn sampler: adaptively setting path lengths in Hamiltonian Monte Carlo *J. Mach. Learn. Res.* **15** 1593–623
[28] Jørgensen J S *et al* 2021 Core Imaging Library—part I: a versatile Python framework for tomographic imaging *Phil. Trans. R. Soc.* A **379** 20200192
[29] Kaipio J and Somersalo E 2005 *Statistical and Computational Inverse Problems* (Springer)
[30] Kumar R, Carrolland C, Hartikainenand A and Martin O 2019 ArviZ a unified library for exploratory analysis of Bayesian models in Python *J. Open Source Softw.* **4** 1143
[31] Laumont R, Bortoli V D, Almansa A, Delon J, Durmus A and Pereyra M 2022 Bayesian imaging using plug & play priors: when Langevin meets Tweedie *SIAM J. Imaging Sci.* **15** 701–37
[32] Li S Z 2009 *Markov Random Field Modeling in Image Analysis* (Springer)
[33] Logg A, Mardal K-A and Wells G 2012 *Automated Solution of Differential Equations by the Finite Element Method—The Fenics Book* (Springer)
[34] Luengo D, Martino L, Bugallo M, Elvira V and Särkkä S 2020 A survey of Monte Carlo methods for parameter estimation *EURASIP J. Adv. Signal Process.* **2020** 25
[35] Marelli S and Sudret B 2014 UQLab: a framework for uncertainty quantification in MATLAB *Vulnerability, Uncertainty, and Risk* pp 2554–63
[36] Meaney A, de Moura F S and Siltanen S 2022 Helsinki tomography challenge 2022 open tomographic dataset (HTC 2022) *Zenodo* (https://doi.org/10.5281/zenodo.6984868)
[37] Nemeth C and Fearnhead P 2021 Stochastic gradient Markov chain Monte Carlo *J. Am. Stat. Assoc.* **116** 433–50
[38] Owen A B 2019 Monte Carlo theory, methods and examples (available at: statweb.stanford.edu/~owen/mc/)
[39] Papoutsellis E, Ametova E, Delplancke C, Fardell G, Jørgensen J S, Pasca E, Turner M, Warr R, Lionheart W R B and Withers P J 2021 Core Imaging Library—part II: multichannel reconstruction for dynamic and spectral tomography *Phil. Trans. R. Soc.* A **379** 20200193
[40] Parno M, Davis A and Seelinger L 2021 MUQ: the MIT uncertainty quantification library *J. Open Source Softw.* **6** 3076
[41] Paszke A, Gross S, Chintala S, Chanan G, Yang E, DeVito Z, Lin Z, Desmaison A, Antiga L and Lerer A 2017 Automatic differentiation in PyTorch *31st Conf. on Neural Information Processing Systems (NIPS 2017)* (*Long Beach, CA, USA, 4–9 December 2017*) pp 1–4
[42] Paszke A *et al* 2019 PyTorch: an imperative style, high-performance deep learning library *Advances in Neural Information Processing Systems* vol 32 pp 8024–35 (Curran Associates, Inc.)
[43] Piponi D, Moore D and Dillon J V 2020 Joint distributions for tensorflow probability (arXiv:2001.11819)
[44] Robert C P and Casella G 2004 *Monte Carlo Statistical Methods* 2nd edn (Springer)
[45] Roberts G O and Rosenthal J S 1998 Two convergence properties of hybrid samplers *Ann. Appl. Probab.* **8** 397–407
[46] Roberts G O and Stramer O 2002 Langevin diffusions and Metropolis-Hastings algorithms *Methodol. Comput. Appl. Probab.* **4** 337–57
[47] Roberts G O and Tweedie R L 1996 Exponential convergence of Langevin distributions and their discrete approximations *Bernoulli* **2** 341–63
[48] Roy V 2020 Convergence diagnostics for Markov Chain Monte Carlo *Annu. Rev. Stat. Appl.* **7** 387–412
[49] Rubin D B 1981 Estimation in parallel randomized experiments *J. Educ. Stat.* **6** 377–401
[50] Rue H and Held L 2005 *Gaussian Markov Random Fields. Theory and Applications* (Chapman & Hall/CRC)
[51] Salvatier J, Wiecki T V and Fonnesbeck C 2016 Probabilistic programming in Python using PyMC3 *PeerJ Comput. Sci.* **2** e55
[52] Stuart A M 2010 Inverse problems: a Bayesian perspective *Acta Numer.* **19** 451–559
[53] Sun S 2013 A review of deterministic approximate inference techniques for Bayesian machine learning *Neural Comput. Appl.* **23** 2039–50

[54] Suuronen J, Chada N K and Roininen L 2022 Cauchy Markov random field priors for Bayesian inversion *Stat. Comput.* **32** 33
[55] Tarantola A 2005 *Inverse Problem Theory and Methods for Model Parameter Estimation* (SIAM)
[56] Tenorio L 2017 *An Introduction to Data Analysis and Uncertainty Quantification for Inverse Problems* (SIAM)
[57] Tran D, Kucukelbir A, Dieng A B, Rudolph M, Liang D and Blei D M 2016 Edward: a library for probabilistic modeling, inference, and criticism (arXiv:1610.09787)
[58] Uribe F, Bardsley J M, Dong Y, Hansen P C and Riis N A B 2022 A hybrid Gibbs sampler for edge-preserving tomographic reconstruction with uncertain view angles *SIAM/ASA J. Uncertain. Quantification* **10** 1293–320
[59] Vehtari A, Gelman A, Simpson D, Carpenter B and Bürkner P-C 2021 Rank-normalization, folding and localization: an improved $\widehat{R}$ for assessing convergence of MCMC (with discussion) *Bayesian Anal.* **16** 667–718
[60] Zhdanov M 2002 *Geophysical Inverse Theory and Regularization Problems* (Elsevier)